\documentclass{article}
\pdfoutput=1
\usepackage{graphicx}
\usepackage{amsmath}
\usepackage{amssymb}
\usepackage{cite,subcaption,wrapfig}

\begin{document}

\title{On inverse problem with phase retrieval for an inclined line in the parabolic approximation}

\author{R.M. Feshchenko$^{1}\footnote{rusl@sci.lebedev.ru}$, I.A. Artyukov$^{1}$, and A.V. Vinogradov$^{1}$}

\maketitle

$^{1}$P.N. Lebedev Physical Institute of RAS, 53 Leninski Prospect, Moscow, 119991, Russia 

\begin{abstract}
The inverse problem of amplitude reconstruction on an inclined line based on the values of amplitude or its module as recorded on semi-infinite line orthogonal to the beam propagation direction is considered within the framework of 2D parabolic equation. It is demonstrated that this inverse problem, in case when the complex image plane amplitude is known, can be reduced to a singular Cauchy type integral equation. The existence of its solutions requires that certain conditions be met but if a solution exists it is necessary unique. The obtained integral equation is then approximated piece-wisely and the resulting linear algebraic system is solved numerically while applying necessary regularization procedures to enhance the stability of its solutions. Finally, an iterative method of phase retrieval is developed and a set of numerical experiments is performed. 
\end{abstract}

\section{Introduction}

The parabolic wave equation (PWE), which was first introduced by Leontovich and Fock more than 60 years ago \cite{fock1965electromagnetic}, is widely used in computational physics and engineering to describe the propagation of paraxial or quasi-paraxial electromagnetic \cite{levy2000parabolic} as  well as in acoustics \cite{lee1995parabolic} beams in free space as well as in inhomogeneous media. It has been successfully applied for solution of complex problems arising in laser physics \cite{sodha1976v}, electromagnetic radiation propagation \cite{levy2000parabolic}, underwater acoustics \cite{lee1995parabolic, spivack1994coherent}, X-ray optics \cite{sakdinawat2010nanoscale}, microscopy and lenseless imaging \cite{chapman2010coherent,thibault2010x}.

One of the underutilized mathematical properties of the PWE is a possibility to express the field amplitude in a part of free space through the initial values of the amplitude specified on an inclined line or plane (depending on the dimension of the problem) or even on arbitrary curve in the 2D space \cite{feshchenko2018propagation}, which is not necessary orthogonal to the beam propagation direction (see \cite{artyukov2014optical, artyukov2014coherent,artyukov2016x} for further details). Another remarkable property of the PWE is its reversibility. For instance, it is always possible to find the initial values of field amplitude on an infinite line or plane orthogonal to the beam's propagation direction based on the values of complex amplitude on a different infinite orthogonal line or plane with the phase retrieval being possible as well. 

On the other hand the inverse problem in case the field amplitude is sought on an inclined line or plane based on its values or values of its module (in case of a phase retrieval problem) on some orthogonal semi-line or semi-plane, i.e. the inversion of the direct problem mentioned in the previous paragraph (inclined parabolic inverse problem -- IPIP), has not been thoroughly considered yet. Nevertheless a singular integral equation was obtained for the case of 2D PWE in one of our previous works \cite{artyukov2014optical} but no attempt was made to solve it at the time.

It should be noted that new and efficient methods for solving the IPIP will be beneficial for many areas of research but especially for the coherent X-ray imaging \cite{paganin2006coherent,nugent2010coherent}, based on the recently developed powerful, versatile and (quasi)-coherent X-ray sources such as compact X-ray lasers \cite{suckewer2009x,ribic2012status, wang2008phase}, free electron lasers \cite{schmuser2008ultraviolet} and X-ray sources based on high order harmonics generation \cite{popmintchev2012bright}. The coherent X-ray imaging possesses a number of advantages over traditional imaging methods including a possibility of lensless imaging with phase retrieval \cite{roy2011lensless}, diffraction imaging \cite{marathe2010coherent} and a sub-picosecond temporal resolution.

In this paper the singular integral equation from \cite{artyukov2014coherent} for the 2D inverse problem is re-derived and its essential properties such as the existence and uniqueness of its solutions are analyzed. A numerical method to solve the 2D IPIP based on a linear piecewise approximation of the initial amplitude sought is developed, in which the singular integral is approximated numerically and the integral equation is then reduced to a linear algebraic system. Applying an appropriate regularization procedure to this linear system, we demonstrate how the method works for a number of model initial amplitudes. 

As the phase retrieval has become a widely used method in diffractive X-ray, optical and electron imaging \cite{latychevskaia2018iterative,meng2009comparison} we also develop an iterative phase retrieval algorithm for the 2D IPIP, which we then demonstrate on a number of numerical experiments. It allows us to reconstruct the initial phase of field amplitude on an inclined line based on the known module of amplitude in an orthogonal image plane.

\section{Direct problem}

Let us briefly review the solution of the direct problem for the 2D parabolic wave equation\cite{levy2000parabolic} 
\begin{equation}
2ik\frac{\partial u}{\partial z}+\frac{\partial^2 u}{\partial x^2}=0,
\label{1a}
\end{equation}
when initial wave field amplitude $u_0$ is specified on an semi-infinite inclined line as defined below. In equation (\ref{1a}) $k=2\pi/\lambda$ is the wave number, $x$ is the transversal coordinate and $z$ is the longitudinal coordinate along the beam propagation direction. Let us assume that the inclined line is defined by the following equations
\begin{equation}
x+z\tan\theta=0,\quad u_0(z)=u(z\tan\theta,z),\quad z<0
\label{1b}
\end{equation}
where $\theta$ is the angle between this inclined line and axis $z$. It is possible to show that the field amplitude in the domain $x>-z\tan\theta$ can be expressed as
\begin{equation}
u(x,z)=(x+z\tan\theta)\sqrt{\frac{k}{2\pi i}}\int\limits_{-\infty}^z\frac{u_0(\xi)}{(z-\xi)^{3/2}}\exp\left[\frac{ik(x+\xi\tan\theta)^2}{2(z-\xi)}\right]\;d\xi.
\label{1r1}
\end{equation}
Changing the variable in integral (\ref{1r1}) to $s=-\xi/\cos\theta$ we arrive at the following expression for $u$
\begin{equation}
u(x,z)=(x\cos\theta+z\sin\theta)\sqrt{\frac{k}{2\pi i}}\int\limits^{\infty}_{-z/\cos\theta}\frac{u_0(s)}{(z+s\cos\theta)^{3/2}}\exp\left[\frac{ik(x-s\sin\theta)^2}{2(z+s\cos\theta)}\right]\;ds,
\label{1s}
\end{equation}
where it follows from \eqref{1b} that $u_0(s)=u_0(-z/\cos\theta)$. Equation \eqref{1s} coincides with equation (4) from \cite{artyukov2014coherent}. If $\theta=0$ and $\zeta=-z-s$ equation \eqref{1s} transforms into a well known in the mathematical physics expression \cite{samarski1964}
\begin{equation}
u(x,z)=x\sqrt{\frac{k}{2\pi i}}\int\limits^{0}_{-\infty}\frac{u_0(-z-\zeta)}{(-\zeta)^{3/2}}\exp\left[-\frac{ikx^2}{2\zeta}\right]\;d\zeta
\label{1t}
\end{equation}
for the calculation of field amplitude at the semi-line $x>0$ parallel to axis $x$ based on the known field amplitude at the part of semi-axis from $-\infty$ to $z$.

In the formulas given above it was implicitly assumed that the field amplitude vanishes when $s\to\infty$, which is equivalent of requiring that there are no field sources at the infinity. In other words, this means that the so-called transparent boundary condition must be satisfied as was shown in \cite{feshchenko2018propagation}.

\section{Inverse problem}

\subsection{Integral equation}

Let us now assume that field amplitude $u(x,z)$ is known on the semi-line $x\ge0,\; z=0$, which we will call the image plane. The inverse problem (which we called IPIP above) is to find initial field amplitude $u$ at the semi-infinite inclined line defined by equations (\ref{1b}). To solve this problem expression (\ref{1r1}) should be multiplied by
\begin{equation}
\sqrt{\frac{2k}{-\pi i(z-z'')}}\exp\left[-\frac{ik(x+z''\tan\theta)^2}{2(z-z'')}\right]
\label{2f1}
\end{equation}
and then integrated by $x$ from $-z\tan\theta$ to $+\infty$. Taking into account equality
\begin{equation}
\lim_{\mu\to0}\frac{1}{x+i\mu}=P\frac{1}{x}-i\pi\delta(x),
\label{2g}
\end{equation}
the final result is the following integral equation for $u_0$
\begin{multline}
u_0(z'')-\frac{i}{\pi}\sqrt{z-z''}\:P\hspace{-8pt}\int\limits_{-\infty}^z \frac{u_0(\zeta)}{z''-\zeta}\exp\left[\frac{ik}{2}\tan^2\theta(z''-\zeta)\right]\frac{d\zeta}{\sqrt{z-\zeta}}=\\
\sqrt{\frac{2ki}{\pi(z-z'')}}\int\limits_{z\tan\theta}^{\infty}u(x',z)\exp\left[-\frac{ik(x'-z''\tan\theta)^2}{2(z-z'')}\right]\;dx',
\label{2f}
\end{multline} 
where the integral in the second term on the left side of (\ref{2f}) is the Cauchy principal value integral. Equation (\ref{2f}) was obtained by us earlier (see \cite{artyukov2014coherent}). 

When $\theta=0$ is assumed in equation (\ref{2f}), the following simplified integral equation is obtained
\begin{multline}
u_0(z'')-\frac{i}{\pi}\sqrt{z-z''}\:P\hspace{-8pt}\int\limits_{-\infty}^z \frac{u_0(\zeta)}{z''-\zeta}\frac{d\zeta}{\sqrt{z-\zeta}}=\\
\sqrt{\frac{2ki}{\pi(z-z'')}}\int\limits_0^{\infty}u(x',z)\exp\left[-\frac{ikx'^2}{2(z-z'')}\right]\;dx'=G(z).
\label{2e}
\end{multline}

From now and on we for the sake of simplicity will discuss equation (\ref{2e}), although all results can be equally applied to equation (\ref{2f}) since the latter follows from (\ref{2e}) by a linear transformation of coordinates: $x'=x+z\tan\theta,\; z'=z$. Equation (\ref{2e}) can be further simplified by assuming $z=0$ and introducing new integration variable $\mu=-\zeta$, new independent variable $t=-z''>0$ and new function $v(t)=u_0(-t)/\sqrt{t}$, which we will also call the amplitude. Now equation (\ref{2e}) can be rewritten as 
\begin{equation}
v(t)+\frac{1}{\pi i}\:P\hspace{-8pt}\int\limits_{0}^{\infty} \frac{v(\mu)}{\mu-t}\,d\mu=
\sqrt{\frac{2ki}{\pi}}\frac{1}{t}\int\limits_0^{\infty}u(x',0)\exp\left[-\frac{ikx'^2}{2t}\right]\;dx'=H(t),
\label{2h}
\end{equation} 
which is a singular integral equation of Cauchy type. 

Singular integral equations similar to \eqref{2h} can be found in a variety of physical and engineering problems. Among them are theory of elasticity and fracture mechanics, hydro- and aerodynamics \cite{ladopoulos2013singular}, electrodynamics and wave mechanics \cite{estrada2012singular} and matter and heat transfer problems \cite{chung1999semi}.

From equation (\ref{2h}) it directly follows that in two special cases -- when the field amplitude $\nu(t)$ is ''a priori'' known to be either a real or imaginary function -- it can be solved exactly by applying complex conjugation to it and adding or subtracting the result from the initial equation. In the real case this gives us the following expression for the field amplitude at the semi-line $x=0,\;z<0$
\begin{equation}
v(t)=\frac{1}{2}\left[H(t)+H^*(t)\right]=\operatorname{Re}H(t),
\label{2i}
\end{equation} 
whereas in the imaginary case we have
\begin{equation}
v(t)=\frac{1}{2}\left[H(t)-H^*(t)\right]=i\operatorname{Im}H(t).
\label{2j}
\end{equation} 
However, in the general case of complex field amplitude, equation (\ref{2h}) can not be solved in that way. 

Let us consider variable $\mu$ to be a complex number. Then the sum of the principle value integral and the amplitude value at point $t$ can be combined into one complex contour integral, where the path is taken around the pole situated on the real axis from the below. So, equation (\ref{2h}) will now take the following form
\begin{equation}
\frac{1}{\pi i}\int\limits_{0}^{\infty} \frac{v(\mu)}{\mu-t}\,d\mu=H(t),
\label{2k}
\end{equation} 
which should be solved relative to function $v$ at the real semi-axis $\mu>0$. From the general theory of Cauchy singular equations it is known (see \cite{estrada2012singular}) that any equation like (\ref{2k}) is a special (generate) case. It has a solution if and only if its right hand side $H(t)$ is a sectionally analytical function with a branch cut along the positive semi-axis. In other words, a solution exists if $H(t)$ is an analytical function everywhere except the positive real semi-axis $t>0$, where it should have a cut. In addition, if a solution of (\ref{2k}) exists it is necessary unique \cite{estrada2012singular}. 

\begin{wrapfigure}{r}{0.5\textwidth}
\includegraphics[width=0.48\textwidth]{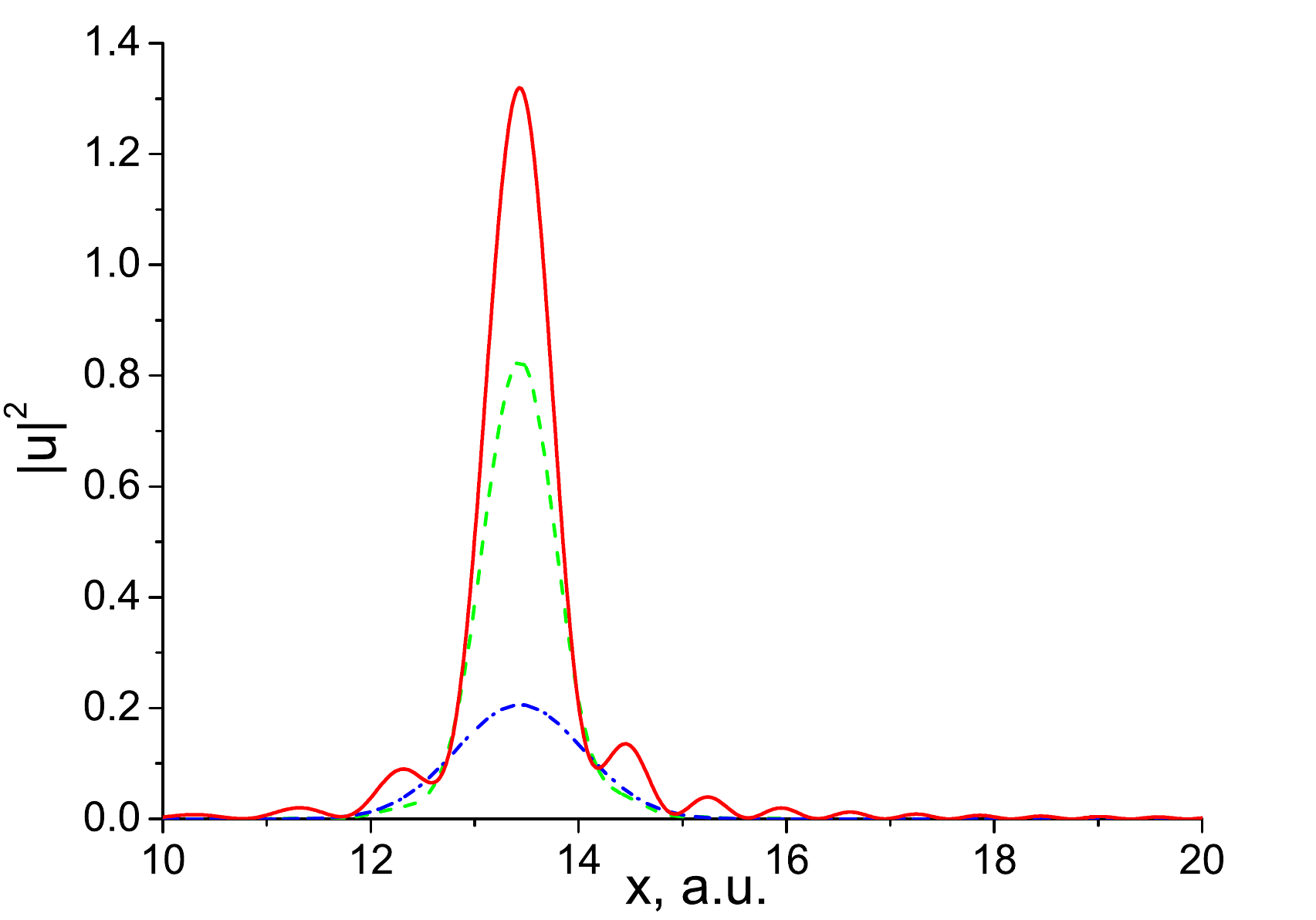}
\caption{The squared module of the image plane field amplitude for the three model fields described by: equation (\ref{3a}) -- blue dash-dot line, equation (\ref{3b}) -- green dash line and equation (\ref{3c}) -- red solid line. The following parameters were used: $z_{min}=90$, $z_{max}=100$, $\tau=0.00398$, $\lambda=0.01$, $\Lambda=1$, $x_{min}=2.8284$,  $x_{max}=28.284$, $N_x=2828$, $N=2514$, $h=0.01$, $L=20$ and $a=4.8$. Arbitrary units were used as units of length and intensity $|u|^2$.}
\label{f1}
\end{wrapfigure}

One can observe that function $H(t)$ as defined in (\ref{2h}) will have a branch cut at the semi-axis $t\ge0$ only when specific conditions are met, for instance, when $u(x,0)=const$. It also may have a singularity at $t=0$. Only in the first case equation (\ref{2k}) will have the necessary unique solution. The absence of solutions for an arbitrary image plane amplitude $u(x,0)$ is a direct consequence of the fact mentioned in the previous section that the amplitude defined by \eqref{1t} must satisfy a transparent boundary condition, and this requirement remain valid for the inverse problem solutions as well. 

From the theory of Cauchy type equations it is also known that in general case the solution, if it exists, is equal to the jump of function $H(t)$ across its branch cut \cite{estrada2012singular}
\begin{equation}
v(t)=H^{+}(t)-H^{-}(t),
\label{2k1}
\end{equation}
where $H^{\pm}(t)$ are the values of $H(t)$ on the upper/lower sides of the cut.

However if only a numerical representation of $u(x,0)$ in a finite interval of $x\in[0,x_0]$ is known, we have to solve equation (\ref{2k}) numerically. This means that real world problems are much more complicated: a numerical solution of (\ref{2k}) may not be unique or it may exist even when the general theory tells us that it should not. 

\subsection{Numerical approximation}
Let us set $z=0$  in equation (\ref{2e}), invert the integration variable as $\zeta\to-\zeta$ and introduce new independent variable $z=-z''>0$. Equation (\ref{2e}) then will take the following form
\begin{equation}
u_0(z)-\frac{1}{\pi i}\sqrt{z}\:P\hspace{-8pt}\int\limits_{0}^{\infty} \frac{u_0(\zeta)}{\zeta-z}\frac{d\zeta}{\sqrt{\zeta}}=\sqrt{\frac{2ki}{\pi z}}\int\limits_{0}^{\infty}u(x',z)\exp\left[-\frac{ik{x'}^2}{2z}\right]\;dx'=G(z).
\label{2k2}
\end{equation} 

Numerical solution of a singular integral equation (e.g. of Cauchy type) can be found using a multitude of methods (see \cite{golberg1990introduction} for a review). They vary with regard to the approximation of the singular integral used. The popular methods include, for instance, the use of Gaussian quadrature rules with various orthogonal polynomials but most often those of Jacobi type. We, however, will use a simple linear piecewise approximation method proposed in \cite{gerasoulis1981method}.

Let us now assume that the field amplitude in the image plane at $z=0$ is known in $N_{x}+1$ points $\{x_i\}$, where $i=0,...N_x$, and that outside interval $[x_0,x_{N_x}]$ the amplitude is zero. Then the function in the right hand side of equation (\ref{2k2}) can be approximated using an appropriate quadrature (for instance, trapezoidal) rule.

\begin{figure}
\begin{subfigure}{.5\textwidth}
 \includegraphics[width=1\linewidth]{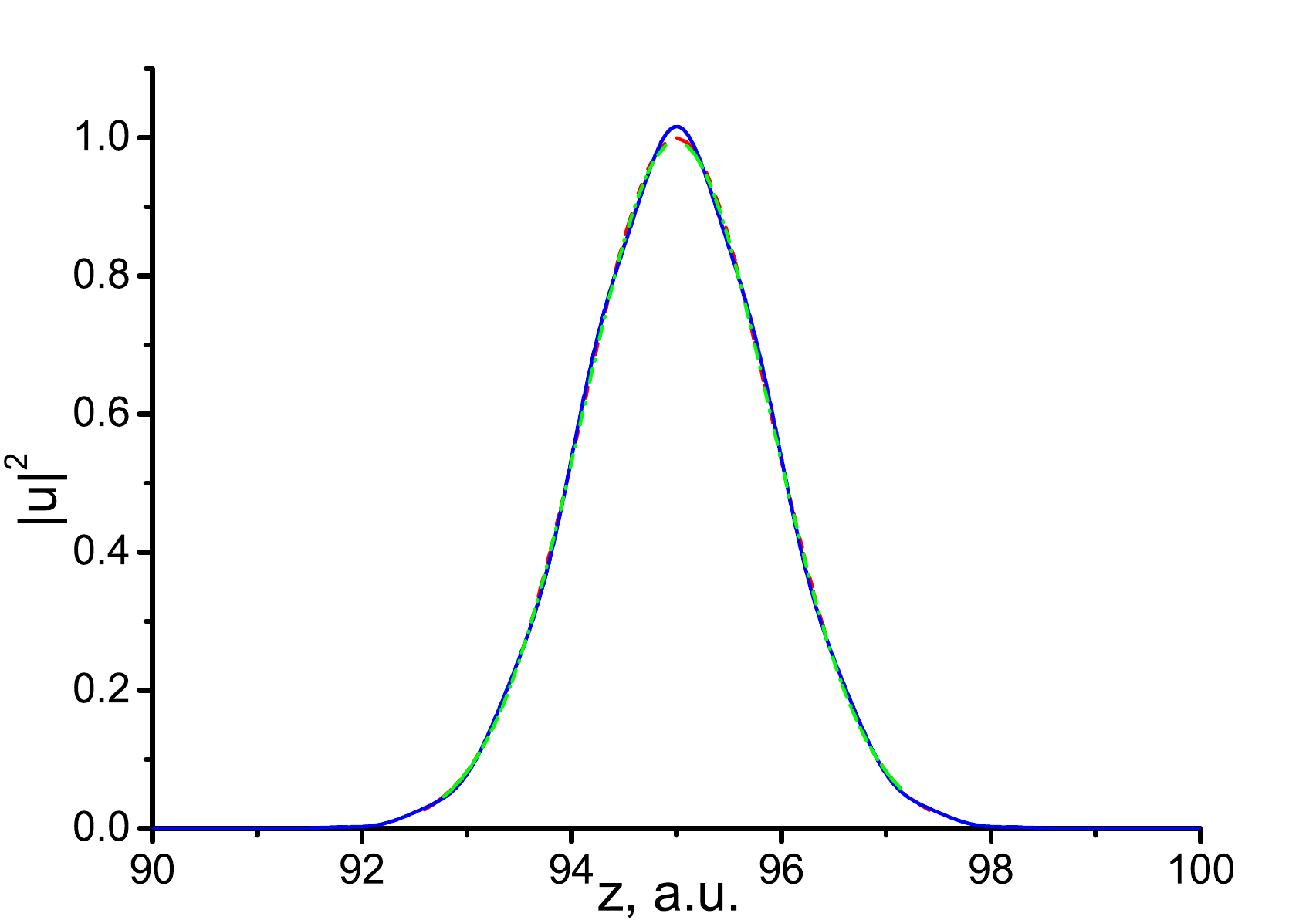}
\end{subfigure}
\begin{subfigure}{.5\textwidth}
 \includegraphics[width=1\linewidth]{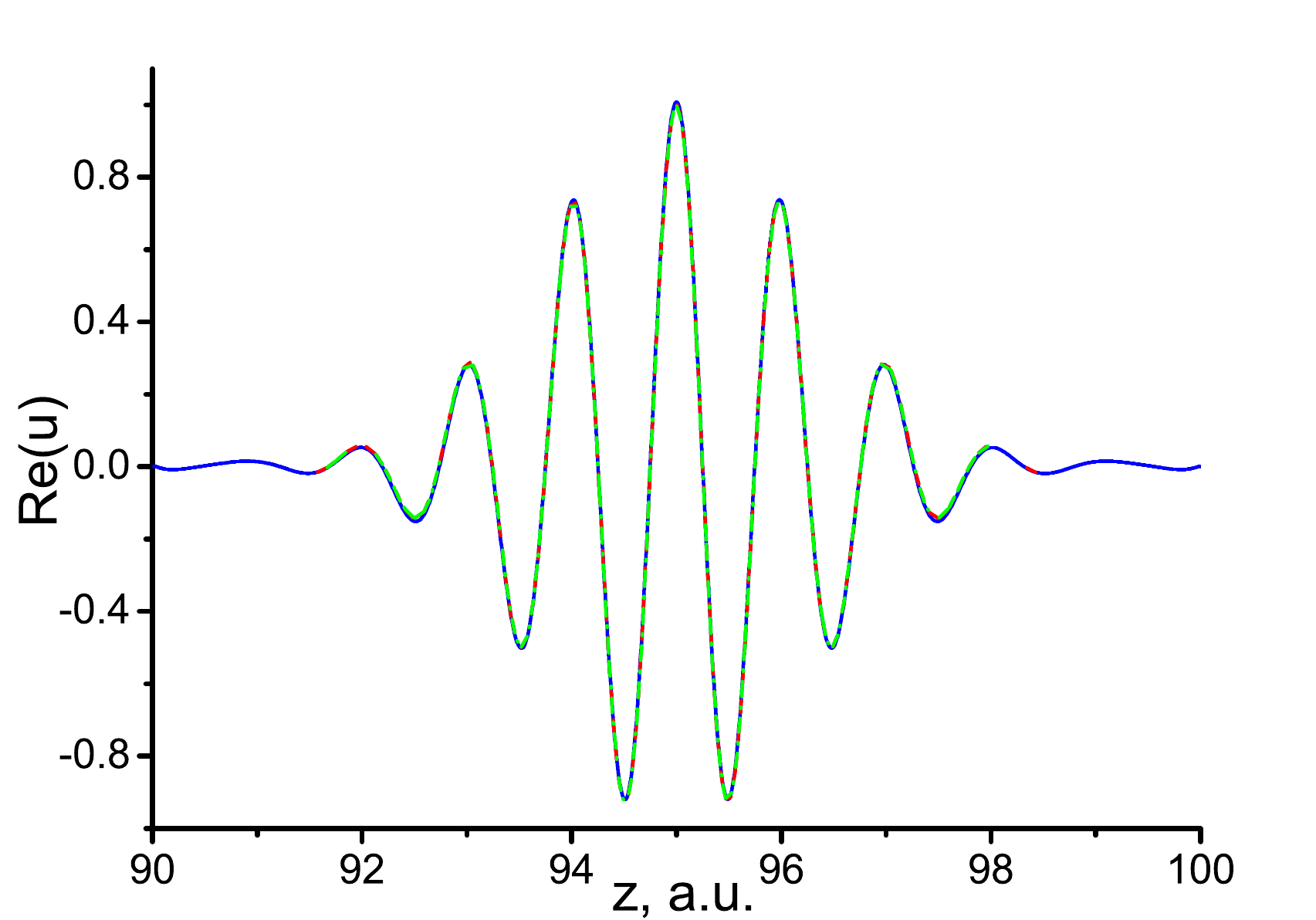}
\end{subfigure}
\begin{subfigure}{.5\textwidth}
 \includegraphics[width=1\linewidth]{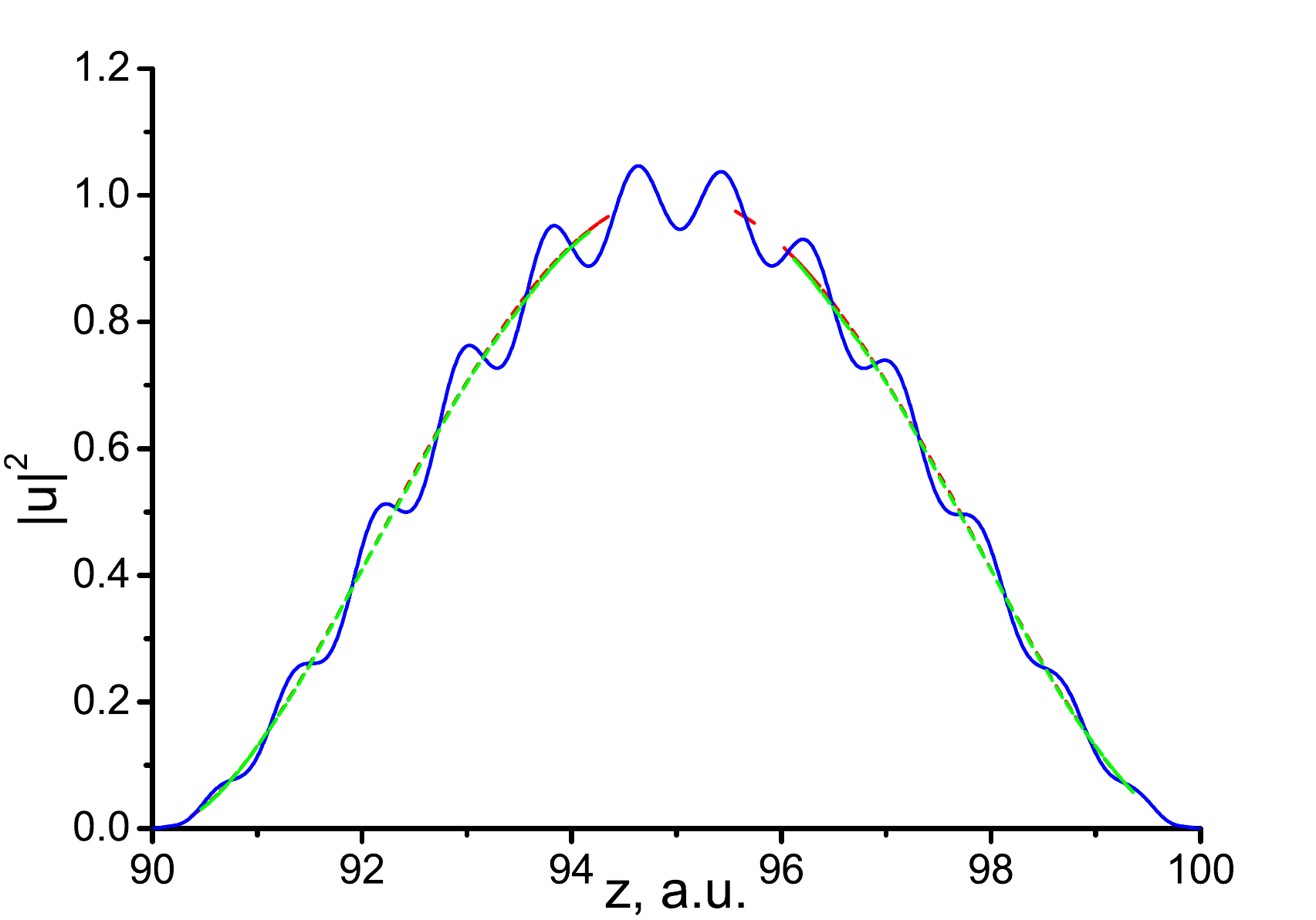}
\end{subfigure}
\begin{subfigure}{.5\textwidth}
 \includegraphics[width=1\linewidth]{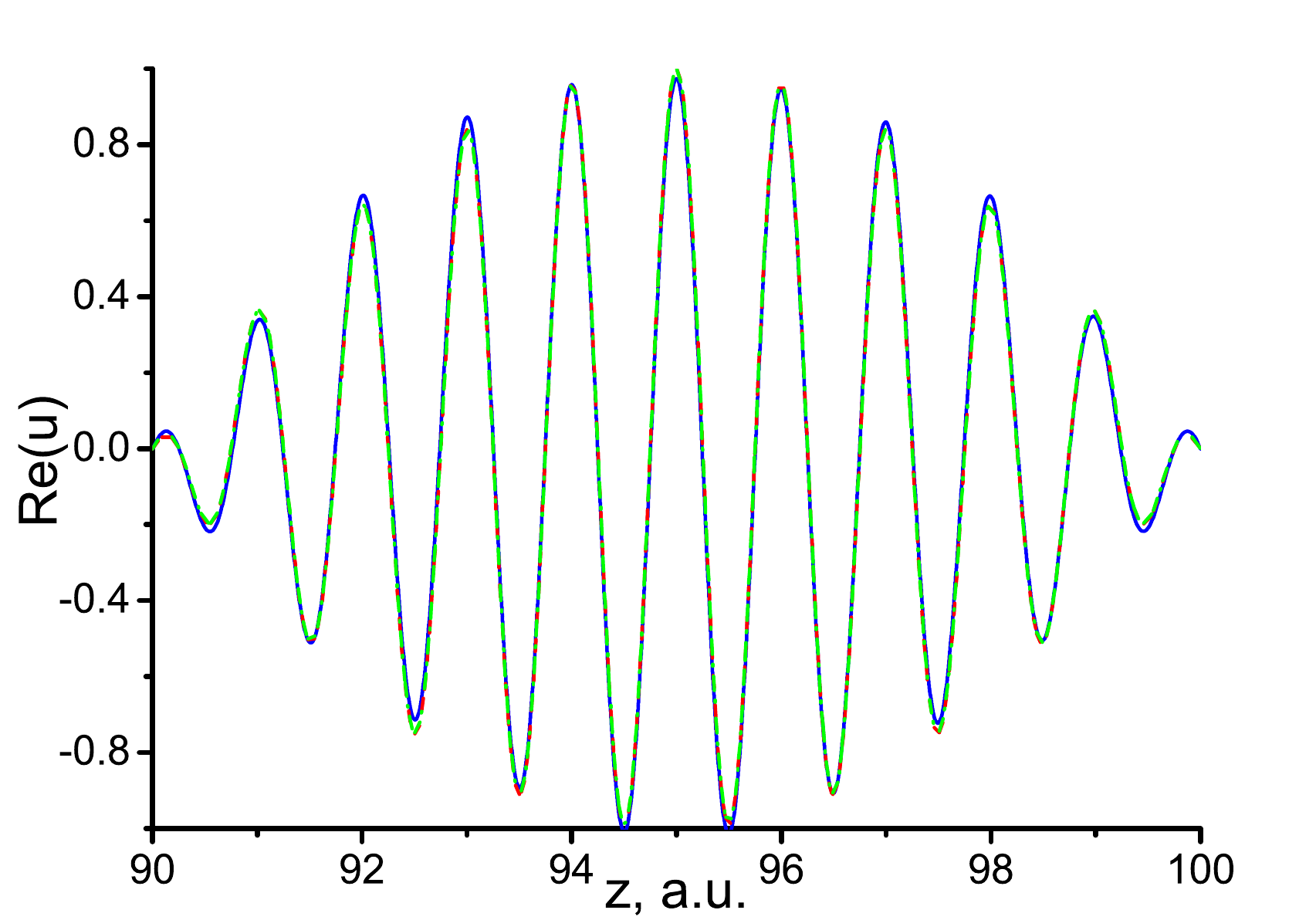}
\end{subfigure}
\begin{subfigure}{.5\textwidth}
 \includegraphics[width=1\linewidth]{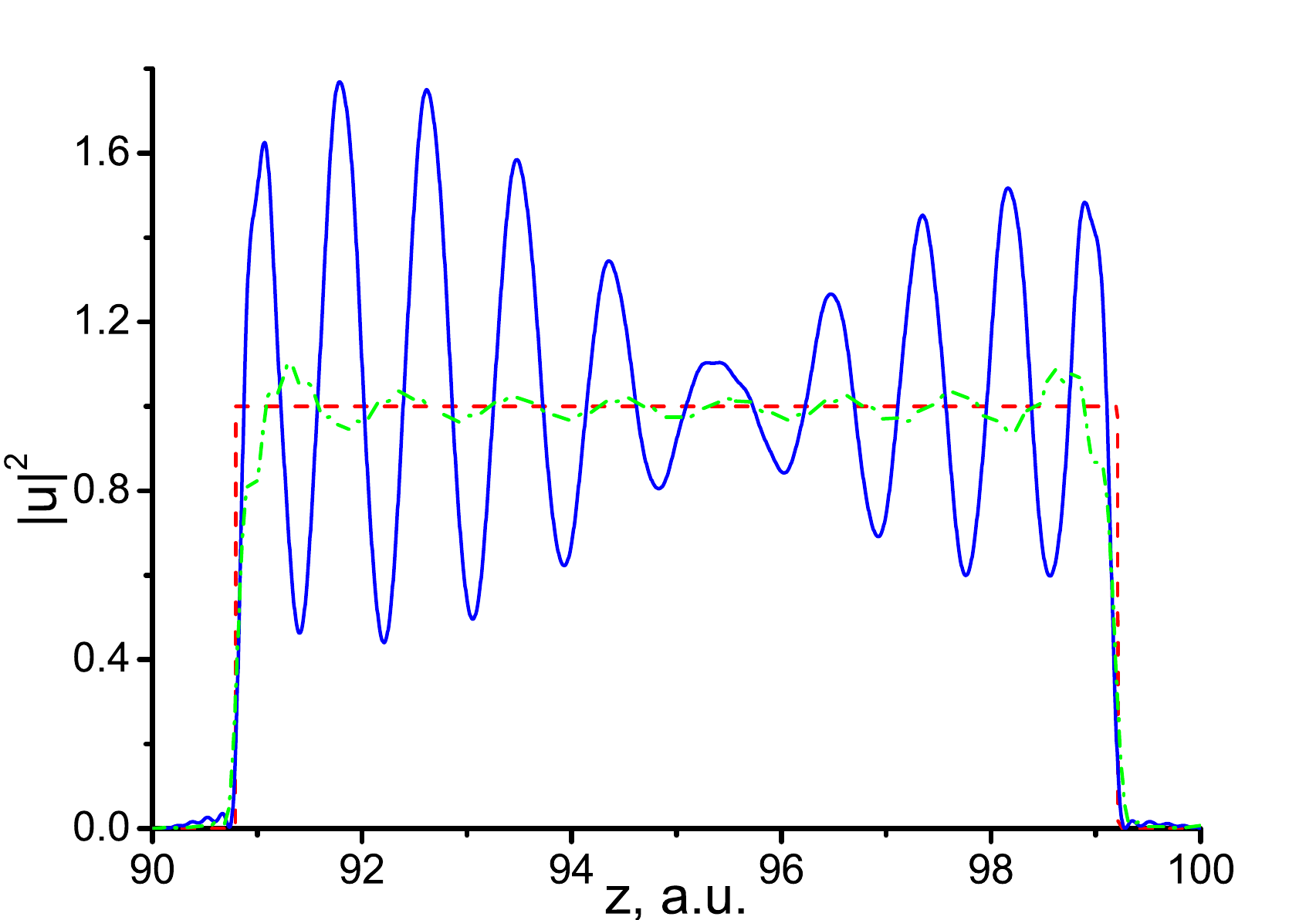}
\end{subfigure}
\begin{subfigure}{.5\textwidth}
 \includegraphics[width=1\linewidth]{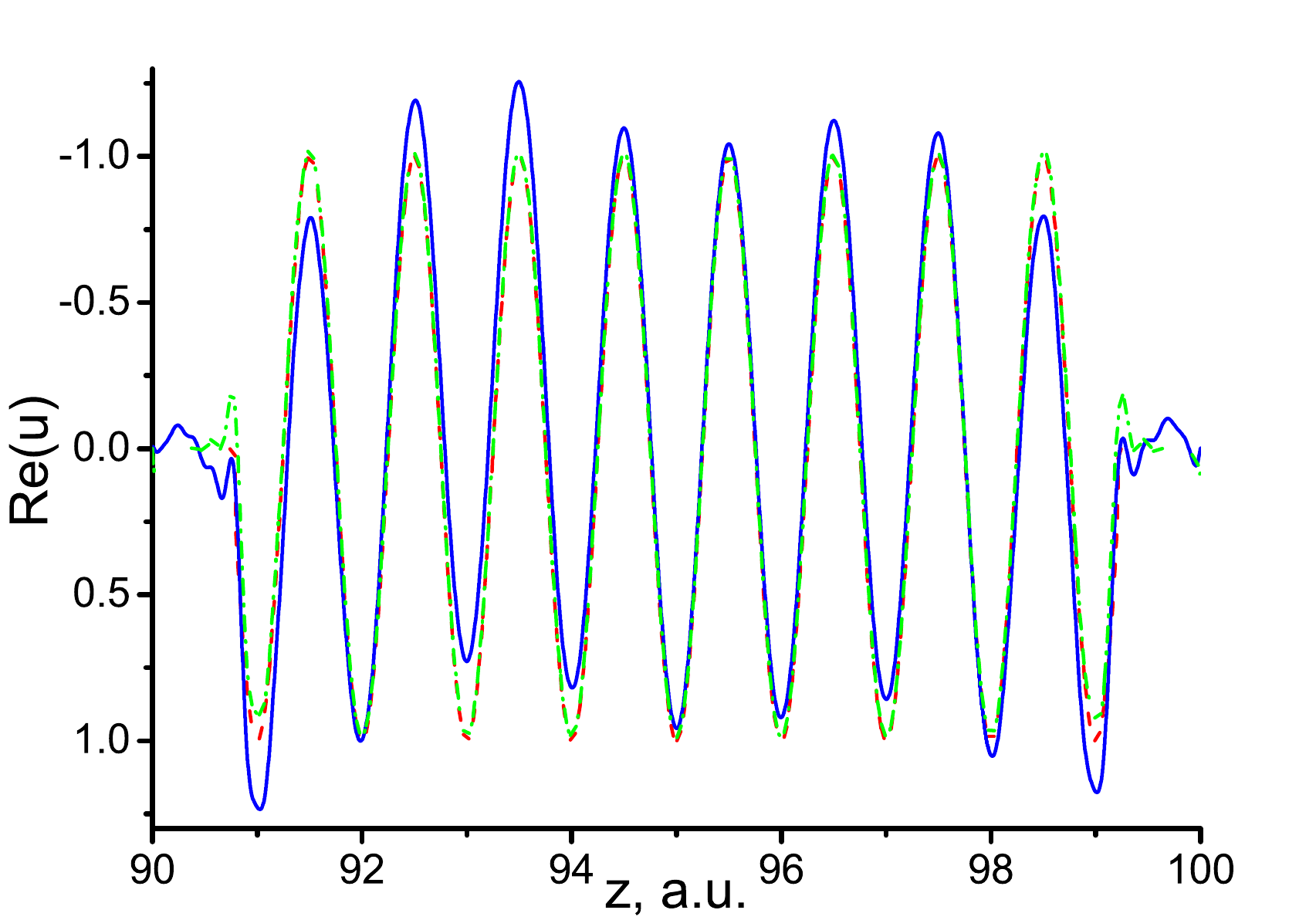}
\end{subfigure}
\caption{The squared module of the initial field amplitude (left panels) and its real part (right panels) for the three model fields described by: equation (\ref{3a}) -- upper panels, equation (\ref{3b}) -- middle panels and equation (\ref{3c}) -- lower panels. The parameters used are the same as in Fig.\ref{f1}. Red dash lines are pre-specified initial amplitudes, blue solid lines are reconstructions using equation (\ref{2p}) and green dot-dash lines are reconstructions using equation (\ref{2p1}) with regularization parameter $\alpha=0.01$.}
\label{f2}
\end{figure}

To approximate the principle value integral in the left hand side of (\ref{2k2}) we can assume that it is {\it ''a priori''} known that initial amplitude $u_0$ deviates from zero only in a finite interval of $z$: $z_{min}\le z \le z_{max}$. This interval is then split into $N$ smaller intervals and the value of the principle-value integral at point $z_n$ is divided into a sum of integrals over intervals $[z_{m},z_{m+1}]$, where $n,m=0,...N$, $z_m=z_{min}+\tau m$ and $\tau=(z_{max}-z_{min})/N$. The result can be written as

\begin{equation}
\:P\hspace{-8pt}\int\limits_{0}^{\infty}\frac{u_0(\zeta)}{\zeta-z_n}\frac{d\zeta}{\sqrt{\zeta}}=\sum_{m=0}^{N-1}\int\limits_{z_m}^{z_{m+1}}\frac{u_0(\zeta)}{\zeta-z_n}\frac{d\zeta}{\sqrt{\zeta}}=\sum_{m=0}^{N-1}I_m(z_n),
\label{2l}
\end{equation}
where
\begin{equation}
I_m(z_n)=\left\{
\begin{array}{rl}
&\:P\hspace{-8pt}\int\limits_{z_m}^{z_{m+1}}\frac{u_0(z_m)-u_0(z_n)+u_{0}'(z_m)(\zeta-z_m)}{\zeta-z_n}\frac{d\zeta}{\sqrt{\zeta}}+\\
&u_0(z_n)\:P\hspace{-8pt}\int\limits_{z_m}^{z_{m+1}}\frac{1}{\zeta-z_n}\frac{d\zeta}{\sqrt{\zeta}},\; m<n,\\
&\:P\hspace{-8pt}\int\limits_{z_m}^{z_{m+1}}\frac{u_0(z_{m+1})-u_0(z_n)+u_{0}'(z_m)(\zeta-z_{m+1})}{\zeta-z_n}\frac{d\zeta}{\sqrt{\zeta}}+\\
&u_0(z_n)\:P\hspace{-8pt}\int\limits_{z_m}^{z_{m+1}}\frac{1}{\zeta-z_n}\frac{d\zeta}{\sqrt{\zeta}},\;m\ge n,
\end{array}\right.
\label{2m1}
\end{equation}
after $u(\zeta)$ is approximated linearly and where
\begin{equation}
u_{0}'(z_m)=\frac{u_0(z_{m+1})-u_0(z_{m})}{z_{m+1}-z_{m}}.
\label{2n}
\end{equation}
In expressions (\ref{2m1}) two separate linear piecewise interpolations of the amplitude $u_0$ were used between points ${z_m}$ and ${z_{m+1}}$ for $m<n$ and $m\le n$. By grouping terms with the same $u_m$ in (\ref{2l}) the following system of linear algebraic equations is finally obtained

\begin{equation}
\left(\mathbf{I}-\frac{\mathbf{M}}{\pi i}\right)\overrightarrow{\mathbf{u_0}}^T=\overrightarrow{\mathbf{g}}^T,
\label{2p}
\end{equation}
where $\mathbf{M}$ is a matrix $(N+1)\times (N+1)$, $\mathbf{I}$  is a unity matrix $(N+1)\times (N+1)$, $\overrightarrow{\mathbf{u_0}}$ is a row vector made of values $u_0(z_n)$ and $\overrightarrow{\mathbf{g}}$ is a row vector made of values $g_n$ of the approximated right hand side $G(z)$ in (\ref{2k2}). The expressions for elements of matrix $\mathbf{M}$ and vector $\overrightarrow{\mathbf{g}}$ can be found in Appendix \ref{ap1}. 

A solution of system (\ref{2p}) will approximate the initial amplitude, although, as it will be seen below, matrix $\mathbf{M}$ has often a high condition number. This high condition number means that some regularization method should be used. One possibility is to solve the Euler system corresponding to \eqref{2p} instead with an additional small regularizing term \cite{liu2018geophysical}. 

Let's introduce the following operators
\begin{align}
\mathbf{A}=&\mathbf{I}-\frac{\mathbf{M}}{\pi i},\label{2t1}\\
\mathbf{A}^+=&\mathbf{I}+\frac{\mathbf{M}^T}{\pi i},\label{2t2}
\end{align}
where matrix $\mathbf{A}^+$ is the Hermitian conjugate matrix with respect to $\mathbf{A}$. 

Equation \eqref{2p} can be transformed into the corresponding Euler system of equations with multiplication of both sides by operator $\mathbf{A}^+$ and then by adding a new constant real term to the resulting equation to improve the stability of its solutions. Finally, equations \eqref{2p} become
\begin{equation}
(\mathbf{A}^+\mathbf{A})\overrightarrow{\mathbf{u_0}}^T+\alpha\overrightarrow{\mathbf{u_0}}^T=\mathbf{A}^+\overrightarrow{\mathbf{g'}}^T,
\label{2p1}
\end{equation}
where $\alpha\to0$ is a real regularization parameter, which value should be chosen to minimize artifacts of inverse problem solutions. One can observe that matrix $\mathbf{A}^+\mathbf{A}$ is obviously Hermitian and the regularization term does not change this fact. The conditional number of the matrix of system \eqref{2p1} should be much lower than that of \eqref{2p}.

\begin{figure}
\begin{subfigure}{.5\textwidth}
 \includegraphics[width=1\linewidth]{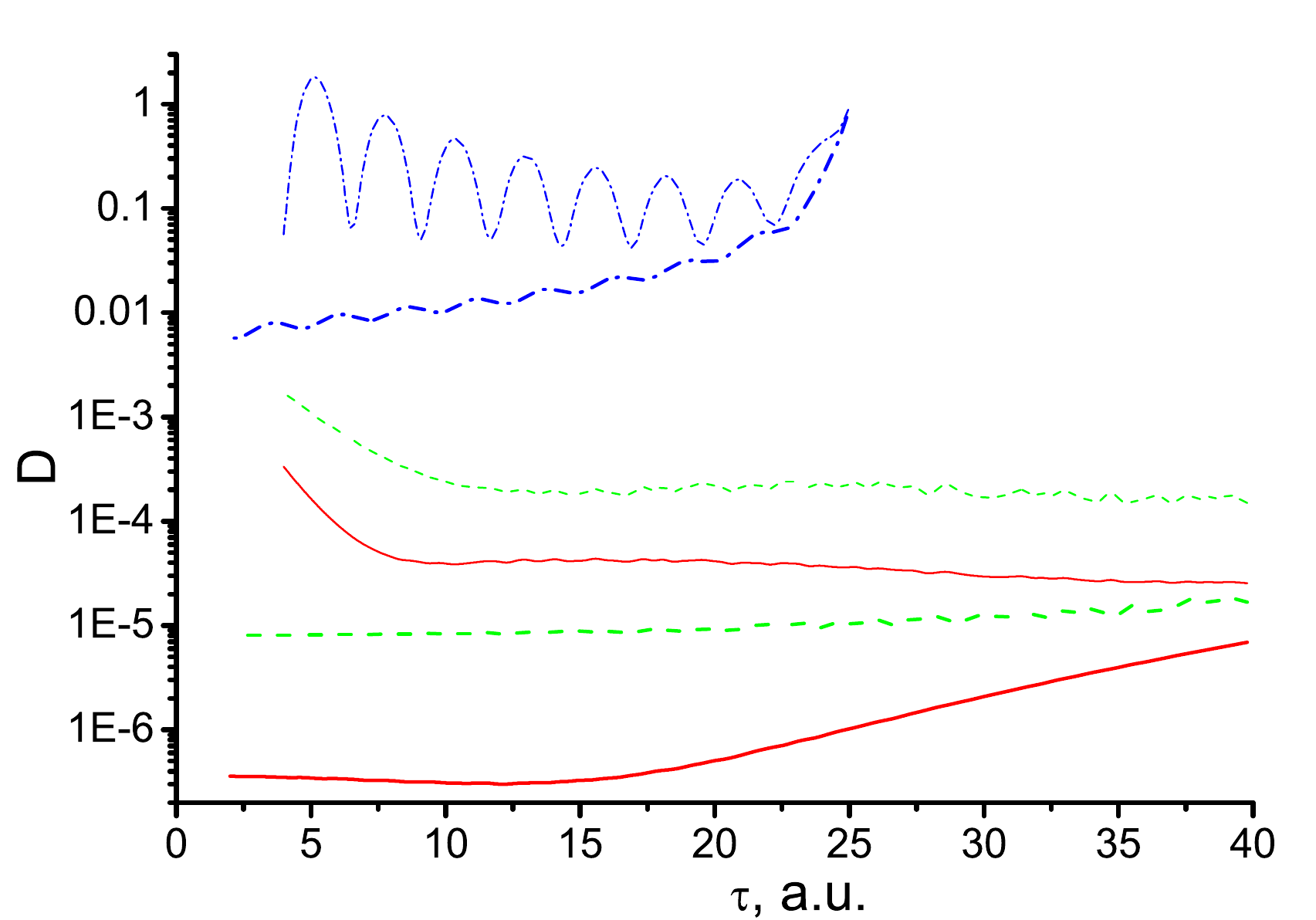}
\end{subfigure}
\begin{subfigure}{.5\textwidth}
 \includegraphics[width=1\linewidth]{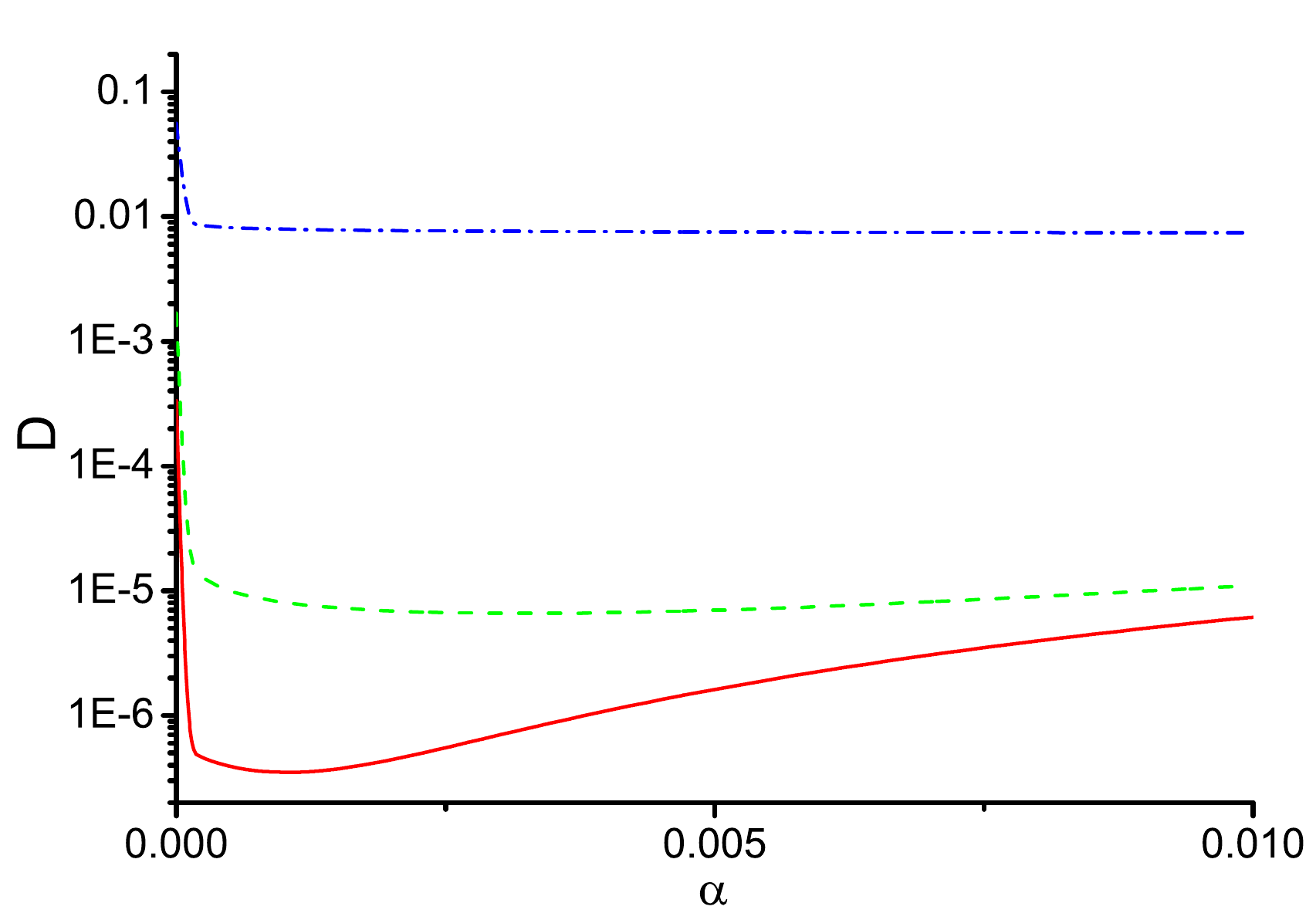}
\end{subfigure}
\caption{Squared difference functions $D(\tau,\alpha)$ for solutions of systems \eqref{2p} and \eqref{2p1} for the Gaussian (solid red curves), parabolic (dash green curves) and step-like (dash-dot blue curves) models as function of the longitudinal step $\tau$ (left) and regularization parameter $\alpha$ (right). The curves for systems \eqref{2p} and \eqref{2p1} are shown by thin and thick lines, respectively. In case of system \eqref{2p1} in the left figure it was also assumed that $\alpha=0.001$. Other parameters are the same in Fig.\ref{f1}. Arbitrary units are used as units of length.}
\label{f3}
\end{figure}

\subsection{Phase retrieval}
The phase retrieval in inverse field propagation problems is commonly achieved using an iterative algorithm, where the known module of the field amplitude in the image plane is supplanted with some specially chosen phase and the resulting complex amplitude is propagated back to the initial plane (see \cite{marchesini2007invited} for a review). The resulting tentative initial amplitude, which is supposed to have a compact support, is then truncated in space by a pre-specified window operator and is propagated forward to the image plane. At this stage the values of module of the amplitude are again replaced with known values and the cycle repeats itself until the sequence of tentative amplitudes converges to a solution. The best known phase retrieval method is error reduction or Gerchberg and Saxton \cite{marchesini2007invited} algorithm, which will be used here for the IPIP solution with phase retrieval.

Let $\hat P$ be the operator (defined in \eqref{1t}) propagating the amplitude from the initial plane to the image plane and $\hat P^{-1}$ be the inverse of it (defined by equation \eqref{2e}). Let function $a(x)$, $x>0$ be a positive real function specifying the module of amplitude in the image plane i.e. $a(x)=|u(x,z=0)|$. Thus, the amplitude in the image plane can be written as 
\begin{equation}
u(x,0)=a(x)e^{i\phi_0(x)},\label{2u}
\end{equation}
where $\phi_0(x)$ is a specially chosen initial phase. Now the iterative sequence of initial amplitudes $\{u^n_0(z)\}$ will take the following form 
\begin{align}
u^0_0(z)=&\hat P^{-1}\left(a(x)e^{i\phi_0(x)}\right),\label{2v1}\\
u^{n+1}_0(z)=&\hat P^{-1}\left(a(x)\frac{\hat P(\hat w(z)u^{n}_0(z))}{|\hat P(\hat w(z)u^{n}_0(z))|}\right),\label{2v2}
\end{align}
where $\hat w(z)$ is the window function satisfying the following conditions $\hat w(0)=0$ and $\hat w(z\to-\infty)=0$. The convergence of sequence \eqref{2v1}--\eqref{2v2} to a solution can be estimated using the following metrics
\begin{equation}
m(n)=\sup_z\frac{|u^{n+1}_0(z)-u^{n}_0(z)|}{|u^{n}_0(z)|}.\label{2v3}
\end{equation}
The iterative sequence $\{u^n_0(z)\}$ is considered converged if condition
\begin{equation}
m(n)<\varepsilon,\label{2v4}
\end{equation}
is satisfied, where $\varepsilon\ll 1$ is a small real number.

The algorithm outlined above was implemented numerically using an approximation of \eqref{1t} by the trapezoidal quadrature rule as operator $\hat P$ (except in case of the Gaussian beam) and regularized system \eqref{2p1} as $\hat P^{-1}$.

\section{Numerical experiments}
\subsection{Complex amplitude reconstruction}
In this section practical applications of equations (\ref{2e}) and \eqref{2p} to IPIP are demonstrated by a number of numerical experiments involving a Gaussian beam, a parabolic beam and a step like initial amplitude $u_0(z)$. 

A Gaussian beam, which is an exact solution of parabolic equation \eqref{1a}, is best suited for testing the existence and stability of IPIP solutions because the amplitude in the image plane can be analytically calculated without a need to numerically solve the direct problem. The following Gaussian beam expression will be used

\begin{multline}
u(x,z)=\frac{1}{\sqrt{1+2i(z-z_c)/kw^2}}\exp(-i\xi x+i\xi^2(z-z_c)/2k)\times\\
\exp\left[-\frac{(x-\xi(z-z_c)/k)^2}{w^2+2i(z-z_c)/k}\right],
\label{3a}
\end{multline}
where $w=\sqrt{2L/k}$ is the Gaussian beam waist radius, $L$ is the Rayleigh length, $z_c$ is the location of the beam waist on axis $z$ and $\xi$ is the spatial frequency of transversal oscillations. On the semi-line at $x=0$ the initial amplitude is $u_0(z)=u(x=0,z)$.

The second beam type, which will be considered here, is the parabolic beam. Its amplitude has the following form on the semi-line at $x=0,\;z<0$
\begin{equation}
u_0(z)=\exp(iKz)\left(l^2-(z-z_c)^2\right)\theta\left(l^2-(z-z_c)^2\right),
\label{3b}
\end{equation}
where $\theta(z)$ is theta (step) function, $l$ is the beam's longitudinal semi-length on $z$ axis, $K=\xi^2/2k=2\pi/\Lambda$ is the spatial frequency of longitudinal oscillations and $\Lambda$ is the wavelength of longitudinal oscillations. It will be assumed here that $z_c=(z_{max}+z_{min})/2$ and $l<(z_{max}-z_{min})/2$. To obtain the field amplitude in the image plane at $z=0$ we will need to solve the direct problem by propagating the field amplitude from $z$ axis using expression (\ref{1t}). This will be done using the trapezoidal integration rule with step $\tau$ defined in the caption of Fig.\ref{f1}. After that the inverse problem can be solved by using the calculated amplitude. The parabolic beam was specifically chosen to vanish at the boundaries of interval $[z_{min},z_{max}]$, which will make the inverse problem solutions more stable.

The last example, which we will consider, is a step-like initial amplitude profile
\begin{equation}
u_0(z)=\exp(iKz)\theta\left(l^2-(z-z_c)^2\right),
\label{3c}
\end{equation}
which as we will show below does not demonstrate a good reversibility. 

\begin{wrapfigure}{r}{0.5\textwidth}
\includegraphics[width=0.48\textwidth]{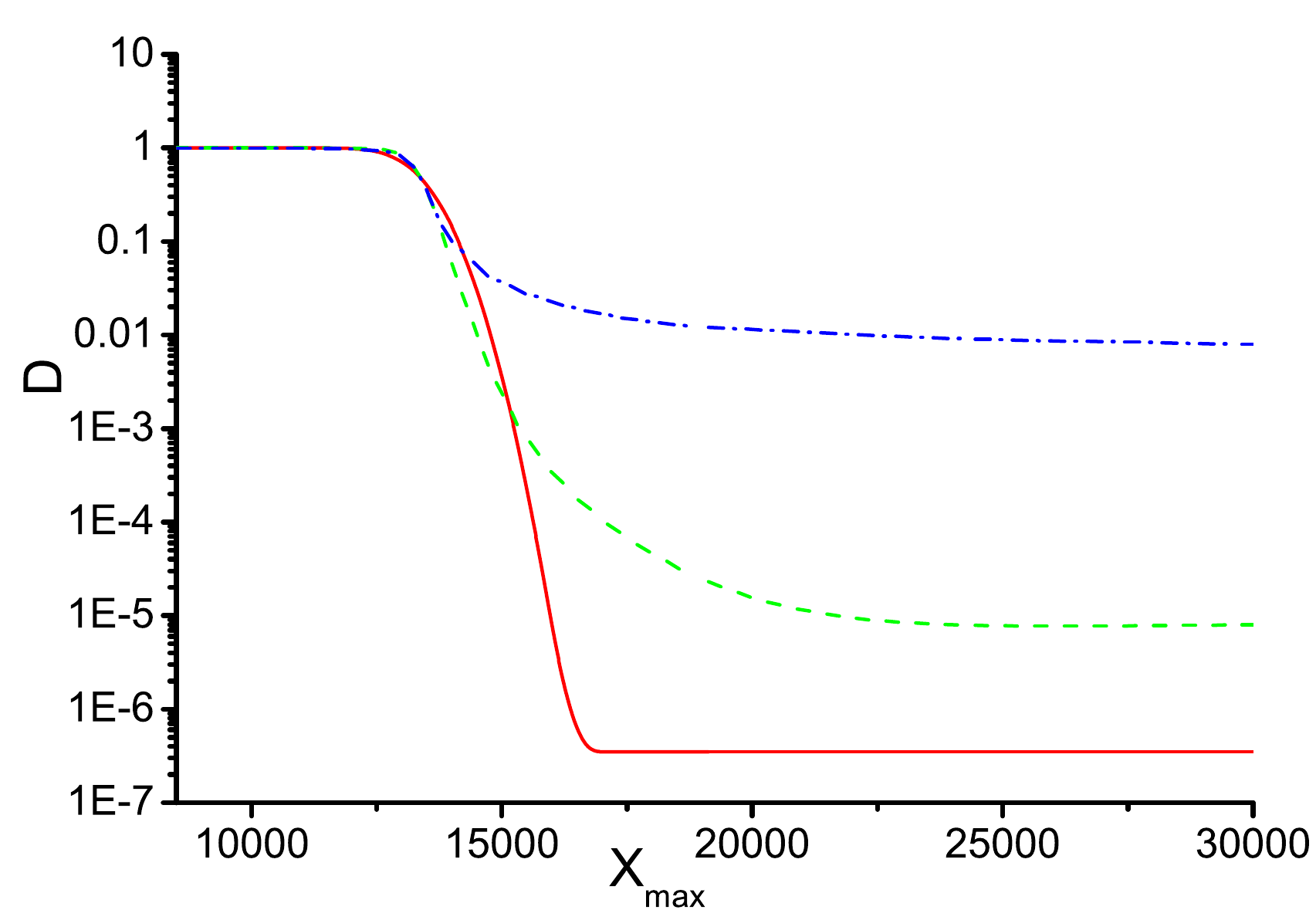}
\caption{Squared difference function $D(x_{max})$ for the Gaussian (solid red curves), parabolic (dash green curves) and step-like (dash-dot blue curves) models as function of the upper limit of integration $x_{max}$ when $\alpha=0.01$ and all other parameters being the same as in Fig.\ref{f1}--\ref{f2}.}
\label{f5}
\end{wrapfigure}

Fig.\ref{f1} shows the amplitudes for these three model fields in the image plane with the numerical simulation parameters specified in its caption. In Fig.\ref{f2} the squared modules of the initial amplitudes and their real parts are shown (as dash red lines) together with the reconstructed values obtained using linear system (\ref{2p}) (as blue solid lines). It can be seen that the reconstructed amplitudes for the Gaussian and parabolic beams approximate the respective initial values relatively well although matrix $\mathbf{M}$ is near singular and its condition number is high -- about $8\cdot 10^{12}$ for all three models (matrix $\mathbf{M}$ does not depend on the model used). However, in case of the step-like initial amplitude there exists a significant difference with the pre-specified initial amplitude particularly near the boarders of computational domain $[z_{min},z_{max}]$. 

The general goodness of numerical solutions is further illustrated by Fig.\ref{f3} where the normalized squared difference function
\begin{equation}
D(\tau, \alpha, X_{max})=\sum_{n=0}^{N}\frac{|u_{0s}(z_n)-u_0(z_n)|^2}{|u_0(z_n)|^2}
\label{3d}
\end{equation}
is plotted as function of longitudinal step $\tau$ (left) and regularization parameter $\alpha$ (right). In \eqref{3d} $u_{0s}$ is the result of solution of either \eqref{2p} or \eqref{2p1}. 

As can be seen in Fig.\ref{f3} the value of difference $D$ for all three models in case when $\alpha=0$ (left, thin lines) reaches a minimum at some small value of $\tau$ and then goes up because the condition number of matrix $\mathbf{M}$ increases for smaller values $\tau$ and the solution becomes worse. In addition, the value of $D$ is always much larger for the step-like model demonstrating that in this case the goodness of solution is worse and that the solution may even be completely incorrect. 

To regularize solutions of system (\ref{2p1}), especially in case of the step-like initial amplitude, we set $\alpha$ to be non-zero. The resulting linear system is not singular any more. The resulting solutions of equation (\ref{2p1}) for the three model initial amplitudes are also shown in Fig.\ref{f2} (as green dash-dot lines) and the respective normalized differences $D$ are shown in Fig.\ref{f3} as function of $\tau$ (left) and $\alpha$ (right) as thick lines. As it can be clearly seen the solution goodness improved noticeably especially for the step-like model.

\begin{figure}
\begin{subfigure}{.5\textwidth}
 \includegraphics[width=1\linewidth]{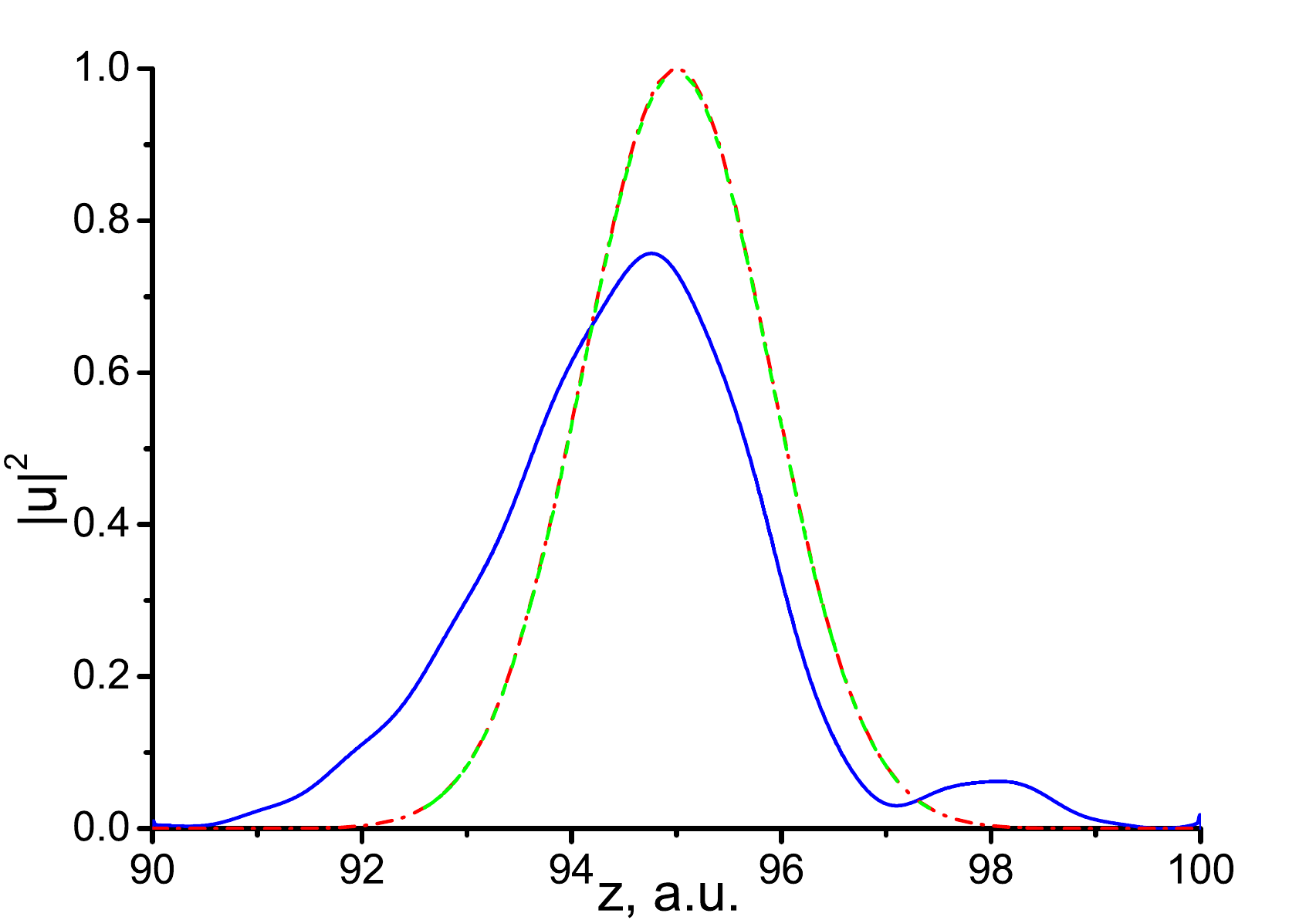}
\end{subfigure}
\begin{subfigure}{.5\textwidth}
 \includegraphics[width=1\linewidth]{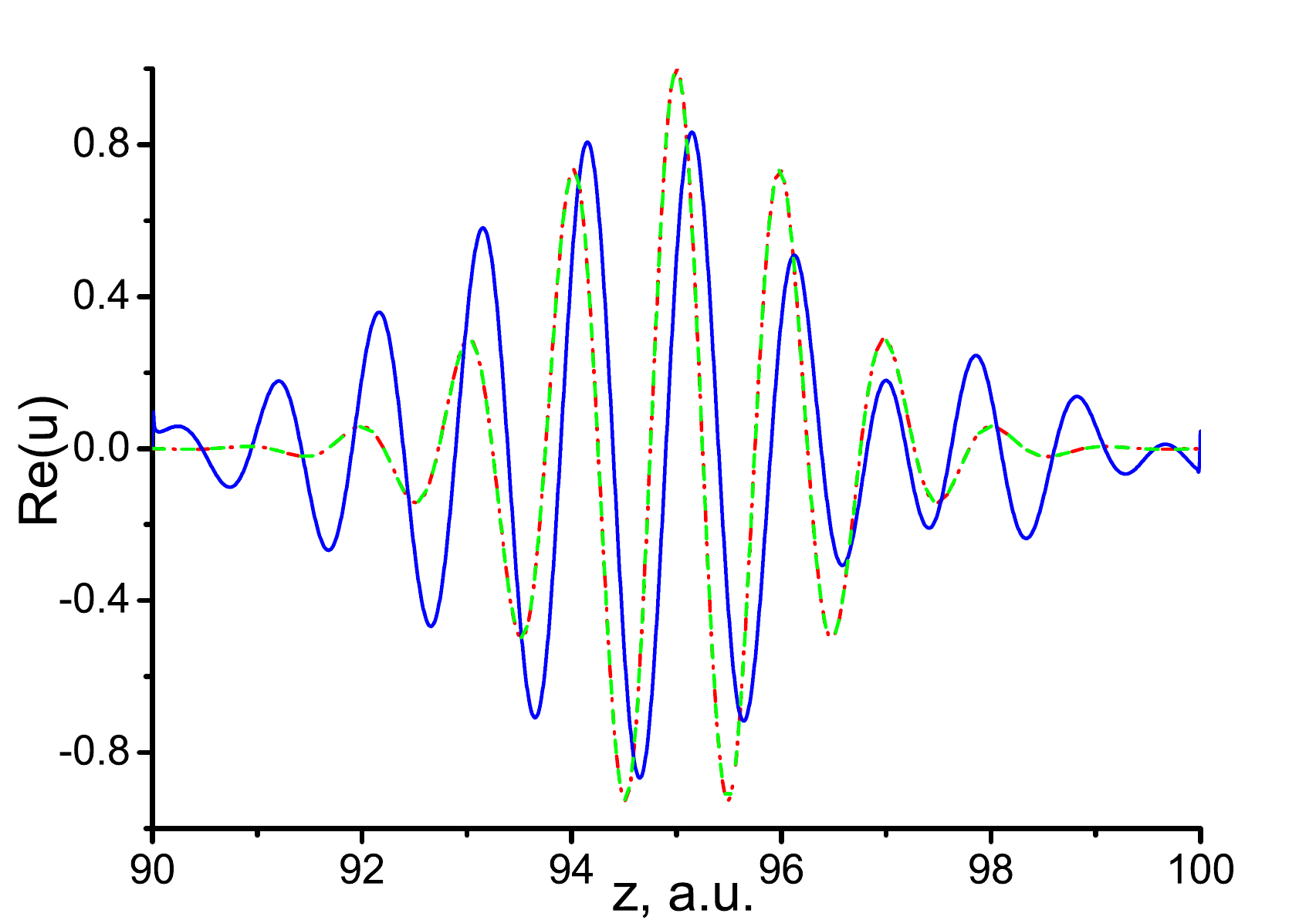}
\end{subfigure}
\begin{subfigure}{.5\textwidth}
 \includegraphics[width=1\linewidth]{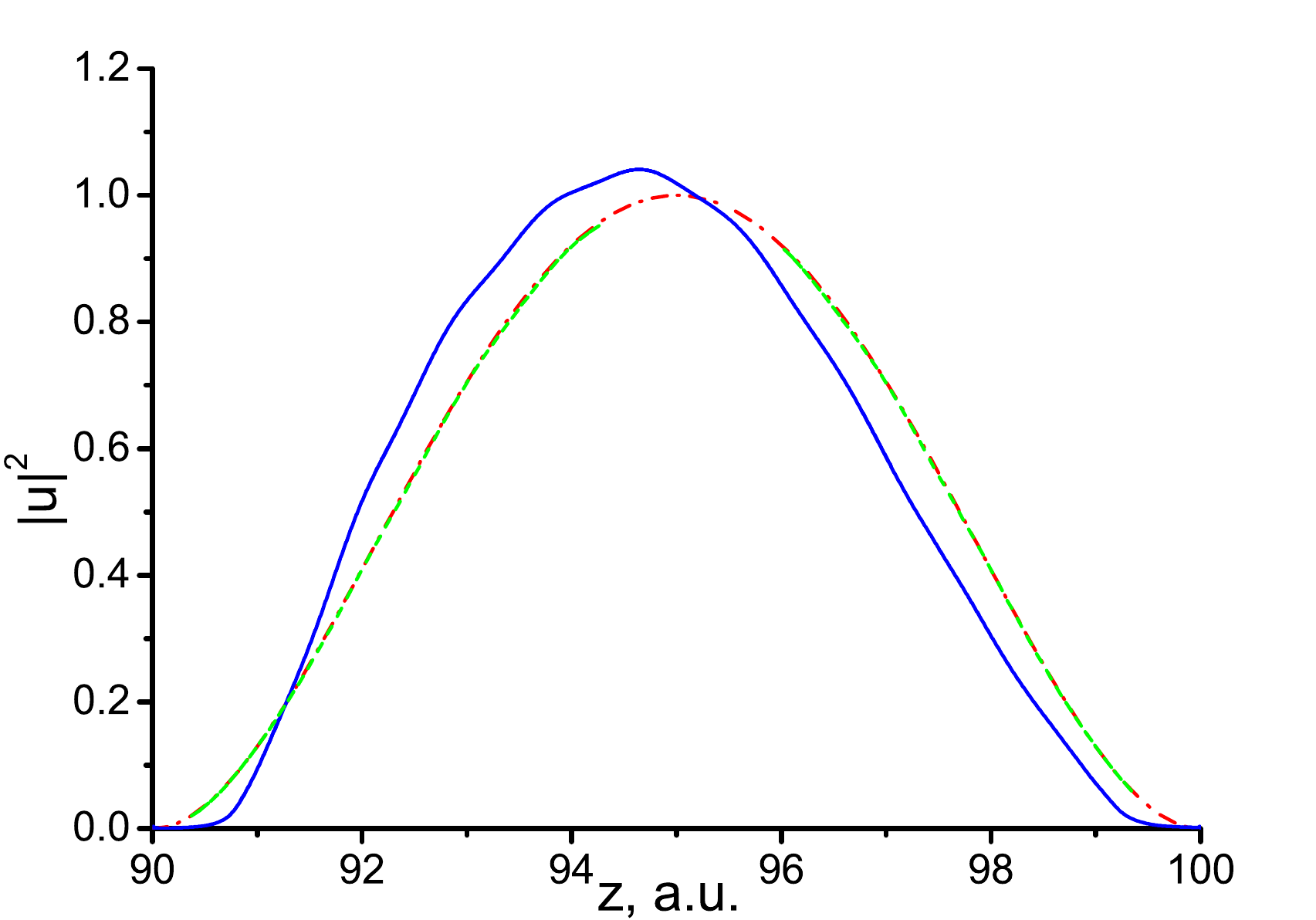}
\end{subfigure}
\begin{subfigure}{.5\textwidth}
 \includegraphics[width=1\linewidth]{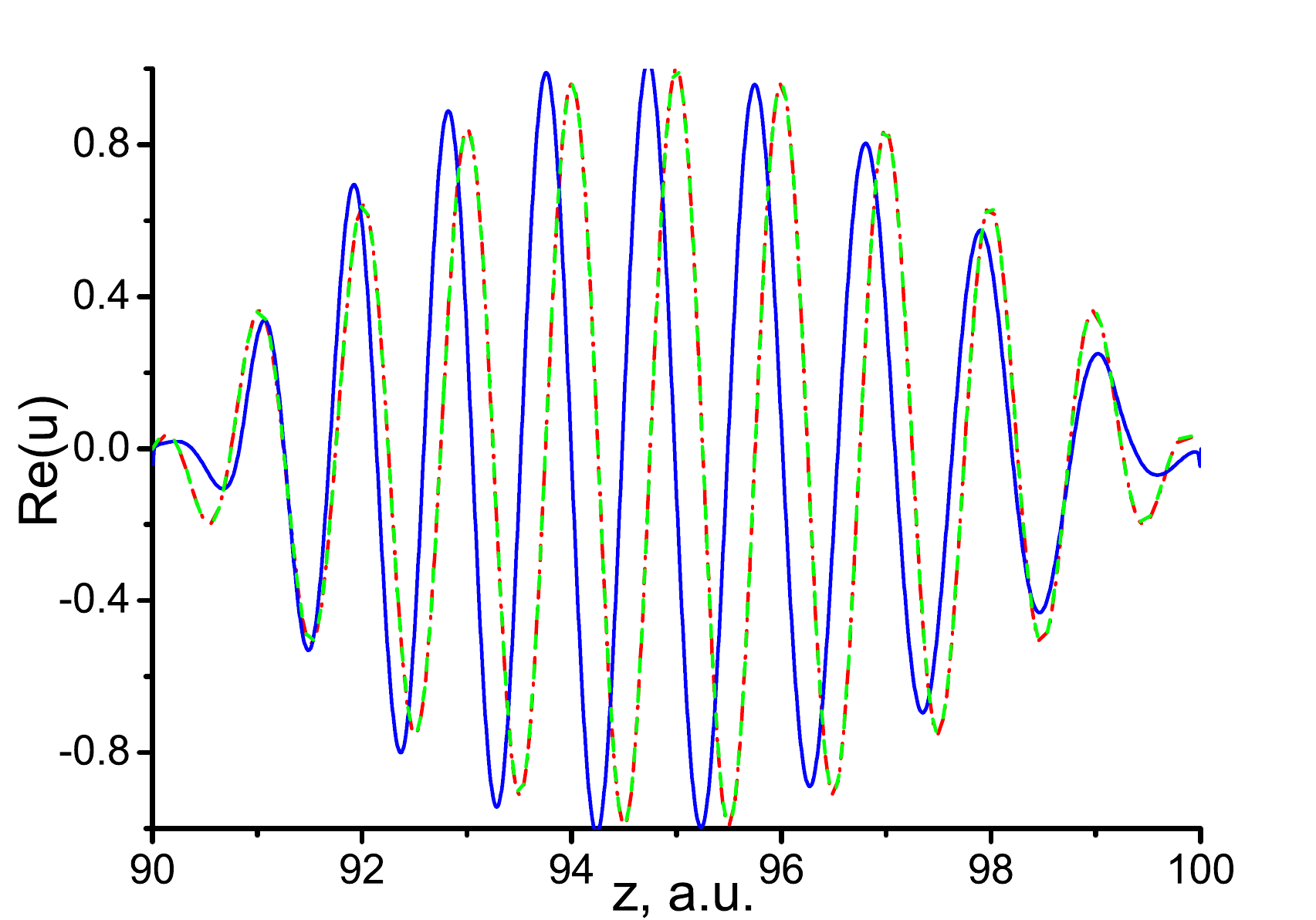}
\end{subfigure}
\begin{subfigure}{.5\textwidth}
 \includegraphics[width=1\linewidth]{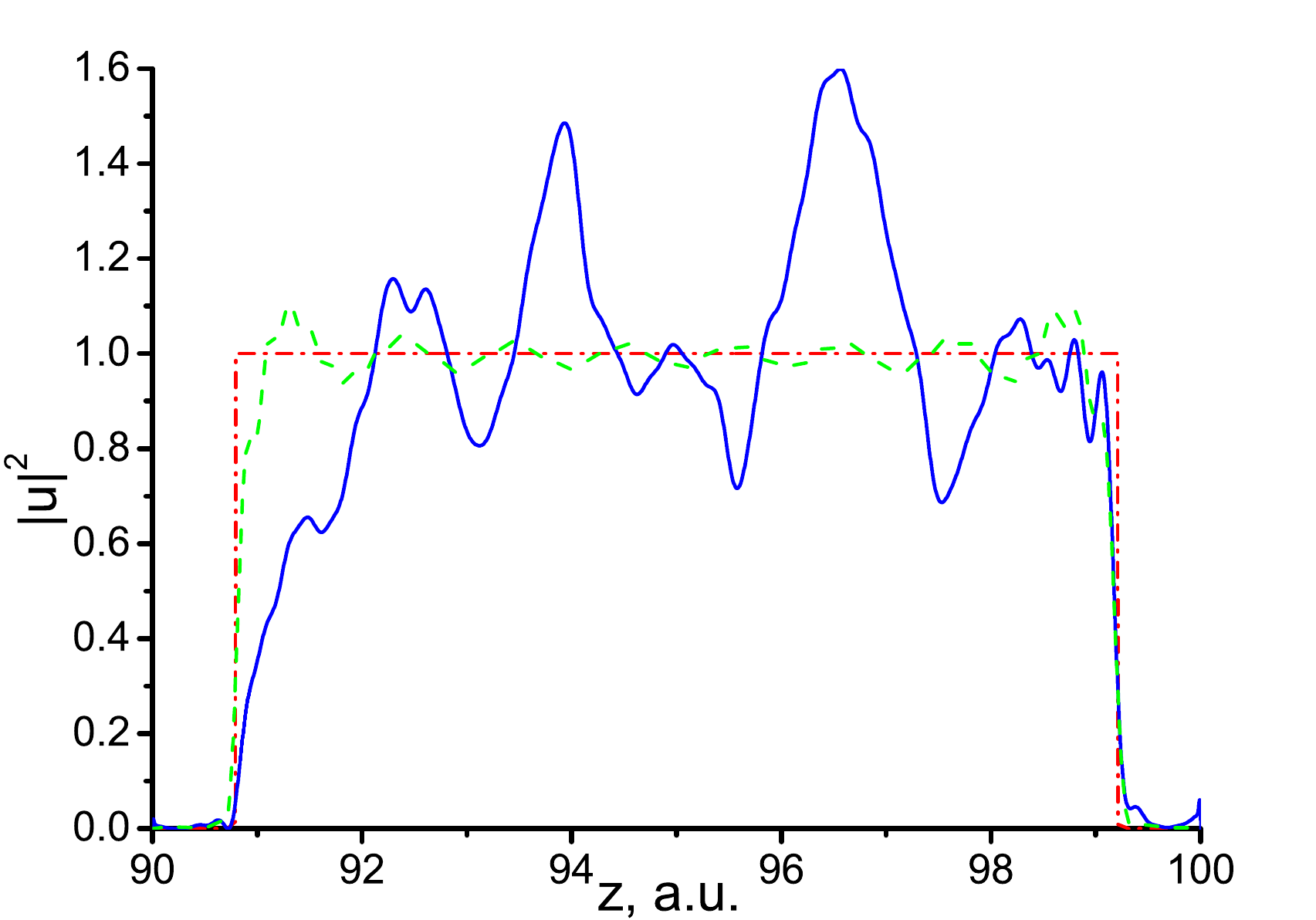}
\end{subfigure}
\begin{subfigure}{.5\textwidth}
 \includegraphics[width=1\linewidth]{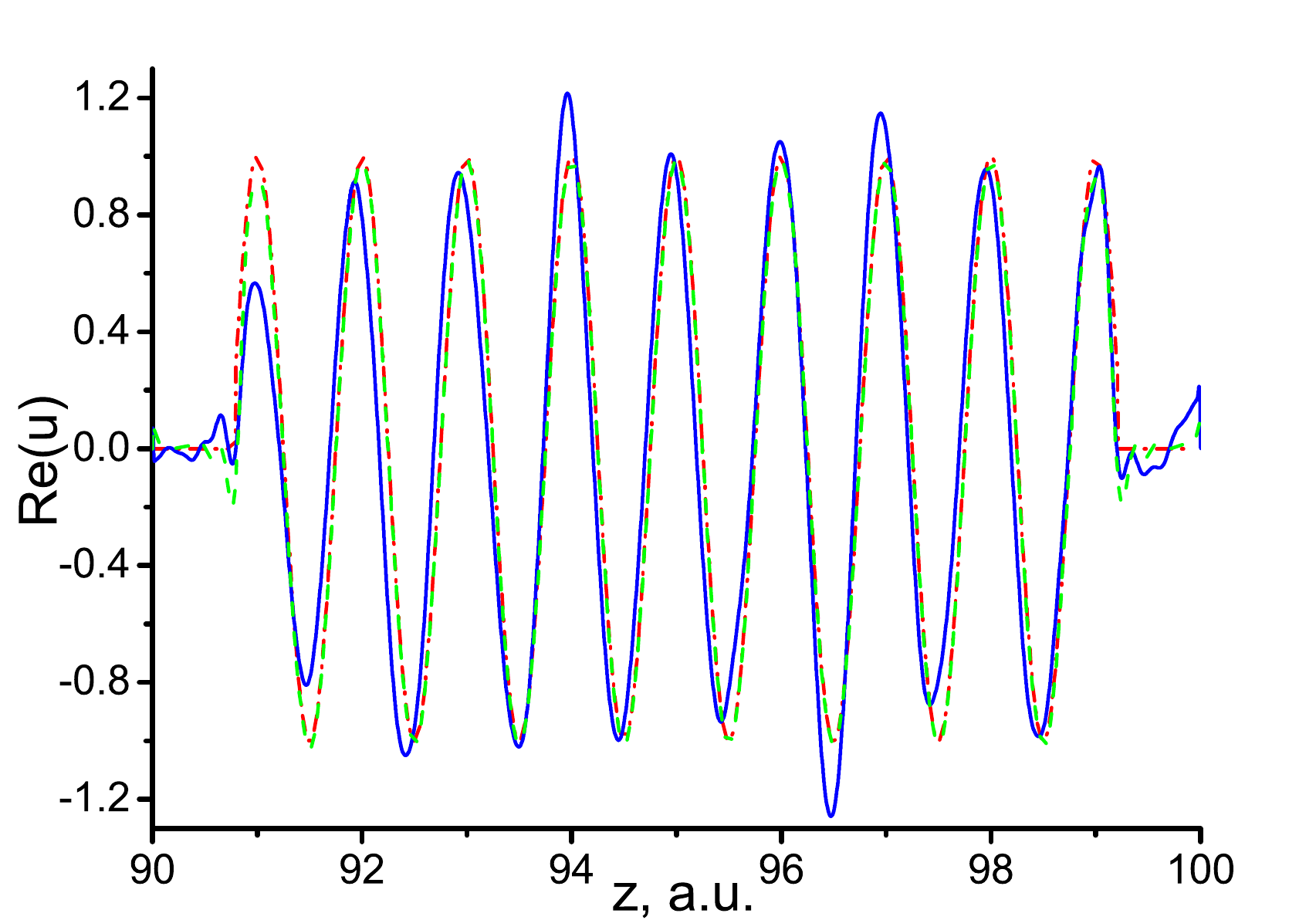}
\end{subfigure}
\caption{The squared module of the initial field amplitude (left panels) and its real part (right panels) for the three model fields described by equations: (\ref{3a}) -- upper panels, equation (\ref{3b}) -- middle panels and equation (\ref{3c}) -- lower panels. The parameters are the same as in Fig.\ref{f1}. For all model initial amplitudes initial phases $\phi_0$ were chosen to be random. The red dash-dot lines are pre-specified initial amplitudes, blue solid lines are the inverse problem solutions with phase retrieval using iterative equations \eqref{2v1}--\eqref{2v2} with $\varepsilon=0.001$ and $\alpha=0.001$. The green dash lines are corresponding solutions of system \eqref{2p1} with $\alpha=0.001$.}
\label{f6}
\end{figure}

It should be noted that in Fig.\ref{f2} and \ref{f3} we chose computational domain $[z_{min}, z_{max}]$ to be wider than the interval where the step-like amplitude deviates from zero, which made the solutions noticeably better. However, from Fig.\ref{f2} and \ref{f3} it is also clear that step-like initial profiles may still not be entirely suitable for the reconstruction. 

From Fig.\ref{f3} (left) it also follows that function $D(\alpha)$ may have a minimum at some value of $\alpha$ when the value of $\tau$ is fixed. As it can be observed in Fig.\ref{f3} (right) there is indeed an optimal value of $\alpha$ for smooth model amplitudes (parabolic and Gaussian), which provides the best approximation. Nevertheless, this is not true for the step-like profile, where the goodness of solutions continues to improve as $\alpha$ increases.

Finally, Fig.\ref{f5} shows the influence of the selection of integration interval $[x_{min},x_{max}]$ on the solution quality. It can be seen that as upper limit $x_{max}$ decreases the deviation $D$ progressively increases until it reaches unity when the upper limit crosses the core part of the image. Thus, from Fig.\ref{f5} it is also clear that the integration interval must include not only the image core but its tail as well, especially for step-like initial amplitudes.

\subsection{Phase retrieval by iterations}
For demonstration of the IPIP with phase retrieval the same three model initial amplitudes were used as for the testing the complex amplitude reversal algorithm in the previous section. The initial phase $\phi_0$ was chosen to be random for all three models. The window operator $\hat w$ in \eqref{2v2} was assumed to have the following step-like form
\begin{equation}
w(z)=\theta\left(l_w^2-(z-z_c)^2\right),
\label{3e}
\end{equation}
where the same notation as in \eqref{3c} is used and the window semi-width is less then the numerical domain being considered i.e. $l_w<l$. The results of the inverse problem solutions with phase retrieval are shown in Fig.\ref{f6} together with solutions of the respective amplitude inverse problems based on system \eqref{2p1}. 

From Fig.\ref{f6} one can see that the quality and veracity of reconstructions is much lower when the phase is unknown, especially in case of the step-like initial model amplitude. In addition, the reconstructed initial amplitudes are somewhat shifted with respect of the initial amplitudes and have additional oscillations superposed on them. The goodness of the phase reconstruction can probably be improved by a different choice of the initial phases because the solution is not unique when the image plane phase is unknown. So, with different initial phases the iterative sequence \eqref{2v1}--\eqref{2v2} may converge to a different solution. 

\section{Summary}

Using the parabolic approximation we solved the 2D IPIP by reversing the previously obtained expression for the field amplitude propagating from an inclined line. We demonstrated that IPIP is reduced to solution of a singular Cauchy type integral equation in a semi-infinite domain. We then analyzed solutions of this equation and found necessary and sufficient conditions for their existence and showed that if a solution exists it is necessary unique. 

We demonstrated that IPIP can be efficiently solved numerically by reducing the obtained singular integral equation to a system of linear algebraic equations and applying the Euler equation method to obtain stable regularized solutions. As in real inverse problems the phase of the image is not usually known, we devised an Gerchberg and Saxton type iterative algorithm for the IPIP with phase retrieval based on the developed amplitude reversal algorithms.

We conducted several numerical experiments using a Gaussian beam, a parabolic-like initial amplitude and a step-like initial amplitude as examples. We demonstrated that for smooth initial field amplitudes of Gaussian and parabolic beams vanishing at the computational domain boarders a good numerical solution of IPIP can be obtained by appropriate selection of the longitudinal step, the image plane integration interval and a regularization parameter $\alpha$ in the corresponding Euler equation. In addition, solving the Euler equation (instead of the initial linear algebraic system) modified by addition of a small term proportional to regularization parameter $\alpha$ allowed us to obtain satisfactory numerical solutions for step-like initial amplitude profiles, for which the quality of non-regularized solutions was initially low. Thus, a careful choice of parameter $\alpha$ is required for better and more stable solutions for both smooth and discontinues initial amplitudes. 

Numerical experiments for the IPIP with phase retrieval for all three above mentioned model amplitudes were also successfully carried out using an iterative phase retrieval algorithm. The obtained numerical solutions resembled the initial amplitudes but had some additional shifts and oscillation as compared to the purely amplitude IPIP solutions. This result shows that the inverse problem solution with phase retrial is feasible for the propagation from inclined objects in the 2D space.   

The developed numerical method for solving IPIP can be used to reconstruct the coherent field amplitude on the surface of an inclined object based on its image recorded within a finite interval in the image plane. In our future work we intend to apply the developed inverse problem solution methods to coherent imaging in the X-ray optics.

\section*{Acknowledgments}

The work was supported by the Basic Research Programme of the Presidium of the Russian Academy of Sciences {\it Actual problems of photonics, probing inhomogeneous mediums and materials}.

\appendix
\section{Matrices and vectors in the linear systems}
\label{ap1}
In this appendix we, for reference, stipulate the elements of matrix $\mathbf{M}$ as well as vector $\overrightarrow{\mathbf{g}}$ used for IPIP solution.

After integration in the formulas (\ref{2m1}) we obtain that
\begin{align}
I_m(z_n)=&2\frac{u_0(z_{m+1})-u_0(z_{m})}{\sqrt{z_{m+1}}+\sqrt{z_{m}}}+\nonumber\\
&\frac{u_0(z_{m})-u_0(z_{n})+(u_0(z_{m+1})-u_0(z_{m}))(z_{n}-z_{m})/\tau}{\tau\sqrt{z_{n}}}\times\nonumber\\
&\ln\left[\frac{(\sqrt{z_{m+1}}-\sqrt{z_{n}})(\sqrt{z_{m}}+\sqrt{z_{n}})}{(\sqrt{z_{m+1}}+\sqrt{z_{n}})(\sqrt{z_{m}}-\sqrt{z_{n}})}\right],\; m<n,\label{2o1}\\
I_m(z_n)&=2\frac{u_0(z_{m+1})-u_0(z_{m})}{\sqrt{z_{m+1}}+\sqrt{z_{m}}}+\nonumber\\
&\frac{u_0(z_{m+1})-u_0(z_{n})+(u_0(z_{m+1})-u_0(z_{m}))(z_{n}-z_{m+1})/\tau}{\tau\sqrt{z_{n}}}\times\nonumber\\
&\ln\left[\frac{(\sqrt{z_{m+1}}-\sqrt{z_{n}})(\sqrt{z_{m}}+\sqrt{z_{n}})}{(\sqrt{z_{m+1}}+\sqrt{z_{n}})(\sqrt{z_{m}}-\sqrt{z_{n}})}\right],\;m\ge n,\label{2o2}
\end{align}

From \eqref{2l} and (\ref{2o1})--(\ref{2o2}) the elements of matrix $\mathbf{M}$ are
\begin{itemize}
\item when $1<m<N-1$ 
\begin{align}
M_{n,N}=&\frac{2}{\sqrt{z_N}+\sqrt{z_{N-1}}}+\frac{n-N+1}{\sqrt{z_n}}\ln\left|\frac{\sqrt{z_N}-\sqrt{z_n}}{\sqrt{z_N}+\sqrt{z_n}}\frac{\sqrt{z_{N-1}}+\sqrt{z_n}}{\sqrt{z_{N-1}}-\sqrt{z_n}}\right|,\label{ap1a1}\\
M_{n,m}=&\frac{4\tau}{(\sqrt{z_m}+\sqrt{z_{m-1}})(\sqrt{z_{m+1}}+\sqrt{z_{m-1}})(\sqrt{z_{m+1}}+\sqrt{z_m})}+\nonumber\\
&\frac{1-n+m}{\sqrt{z_n}}\ln\left|\frac{\sqrt{z_{m+1}}-\sqrt{z_n}}{\sqrt{z_{m+1}}+\sqrt{z_n}}\frac{\sqrt{z_{m-1}}-\sqrt{z_n}}{\sqrt{z_{m-1}}+\sqrt{z_n}}\frac{(\sqrt{z_m}+\sqrt{z_n})^2}{(\sqrt{z_m}-\sqrt{z_n})^2}\right|+\nonumber\\
&\frac{2}{\sqrt{z_n}}\ln\left|\frac{\sqrt{z_m}-\sqrt{z_n}}{\sqrt{z_m}+\sqrt{z_n}}\frac{\sqrt{z_{m-1}}+\sqrt{z_n}}{\sqrt{z_{m-1}}-\sqrt{z_n}}\right|,\quad n+1<m<N\label{ap1a2}\\
M_{n,n+1}=&\frac{4\tau}{(\sqrt{z_{n+1}}+\sqrt{z_n})(\sqrt{z_{n+2}}+\sqrt{z_n})(\sqrt{z_{n+2}}+\sqrt{z_{n+1}})}+\nonumber\\
&\frac{2}{\sqrt{z_n}}\ln\left|\frac{\sqrt{z_{n+2}}-\sqrt{z_n}}{\sqrt{z_{n+2}}+\sqrt{z_n}}\frac{\sqrt{z_{n+1}}+\sqrt{z_n}}{\sqrt{z_{n+1}}-\sqrt{z_n}}\right|,\label{ap1a3}\\
M_{n,n}=&\frac{4\tau}{(\sqrt{z_n}+\sqrt{z_{n-1}})(\sqrt{z_{n+1}}+\sqrt{z_{n-1}})(\sqrt{z_{n+1}}+\sqrt{z_n})}-\nonumber\\
&\frac{1}{\sqrt{z_n}}\ln\left|\frac{\sqrt{z_{n+1}}+\sqrt{z_n}}{\sqrt{z_{n+1}}-\sqrt{z_n}}\frac{\sqrt{z_{n-1}}-\sqrt{z_n}}{\sqrt{z_{n-1}}+\sqrt{z_n}}\right|,\label{ap1a4}\\
M_{n,n-1}=&\frac{4\tau}{(\sqrt{z_{n-1}}+\sqrt{z_{n-2}})(\sqrt{z_n}+\sqrt{z_{n-2}})(\sqrt{z_n}+\sqrt{z_{n-1}})}+\nonumber\\
&\frac{2}{\sqrt{z_n}}\ln\left|\frac{\sqrt{z_{n-1}}-\sqrt{z_n}}{\sqrt{z_{n-1}}+\sqrt{z_n}}\frac{\sqrt{z_{n-2}}+\sqrt{z_n}}{\sqrt{z_{n-2}}-\sqrt{z_n}}\right|,\label{ap1a5}
\end{align}

\begin{align}
M_{n,m}=&\frac{4\tau}{(\sqrt{z_m}+\sqrt{z_{m-1}})(\sqrt{z_{m+1}}+\sqrt{z_{m-1}})(\sqrt{z_{m+1}}+\sqrt{z_m})}+\nonumber\\
&\frac{n-m+1}{\sqrt{z_n}}\ln\left|\frac{\sqrt{z_{m+1}}+\sqrt{z_n}}{\sqrt{z_{m+1}}-\sqrt{z_n}}\frac{\sqrt{z_{m-1}}+\sqrt{z_n}}{\sqrt{z_{m-1}}-\sqrt{z_n}}\frac{(\sqrt{z_m}-\sqrt{z_n})^2}{(\sqrt{z_m}+\sqrt{z_n})^2}\right|+\nonumber\\
&\frac{2}{\sqrt{z_n}}\ln\left|\frac{\sqrt{z_m}-\sqrt{z_n}}{\sqrt{z_m}+\sqrt{z_n}}\frac{\sqrt{z_{m-1}}+\sqrt{z_n}}{\sqrt{z_{m-1}}-\sqrt{z_n}}\right|,\quad 1<m<n-1\label{ap1a6}\\
M_{n,0}=&-\frac{2}{\sqrt{z_1}+\sqrt{z_0}}-\frac{n-1}{\sqrt{z_n}}\ln\left|\frac{\sqrt{z_1}-\sqrt{z_n}}{\sqrt{z_1}+\sqrt{z_n}}\frac{\sqrt{z_0}+\sqrt{z_n}}{\sqrt{z_0}-\sqrt{z_n}}\right|,\label{ap1a7}
\end{align}

\item when $n=N$
\begin{align}
M_{N,N}=&1,\label{ap1b1}\\
M_{N,N-1}=&\frac{4\tau}{(\sqrt{z_{N-1}}+\sqrt{z_{N-2}})(\sqrt{z_N}+\sqrt{z_{N-2}})(\sqrt{z_N}+\sqrt{z_{N-1}})}+\nonumber\\
&\frac{2}{\sqrt{z_N}}\ln\left|\frac{\sqrt{z_{N-1}}-\sqrt{z_N}}{\sqrt{z_{N-1}}+\sqrt{z_N}}\frac{\sqrt{z_{N-2}}+\sqrt{z_N}}{\sqrt{z_{N-2}}-\sqrt{z_N}}\right|.\label{ap1b2}
\end{align}
\item when $n=N-1$
\begin{align}
M_{N-1,N}=&\frac{2}{\sqrt{z_N}+\sqrt{z_{N-1}}},\label{ap1c1}\\
M_{N-1,N-1}=&\frac{4\tau}{(\sqrt{z_{N-1}}+\sqrt{z_{N-2}})(\sqrt{z_N}+\sqrt{z_{N-2}})(\sqrt{z_N}+\sqrt{z_{N-1}})}-\nonumber\\
&\frac{1}{\sqrt{z_{N-1}}}\ln\left|\frac{\sqrt{z_{N-2}}-\sqrt{z_{N-1}}}{\sqrt{z_{N-2}}+\sqrt{z_{N-1}}}\frac{\sqrt{z_N}+\sqrt{z_{N-1}}}{\sqrt{z_N}-\sqrt{z_{N-1}}}\right|.\label{ap1c2}
\end{align}

\item In case $n=N,N-1$ and $m\ne N,N-1$ the elements are determined by general formulas (\ref{ap1a1})--(\ref{ap1a7}). 

\item In the opposite case when $n=0$ the elements are
\begin{align}
M_{0,0}=&-1,\label{ap1d1}\\
M_{0,1}=&\frac{4\tau}{(\sqrt{z_1}+\sqrt{z_0})(\sqrt{z_2}+\sqrt{z_0})(\sqrt{z_2}+\sqrt{z_1})}+\nonumber\\
&\frac{2}{\sqrt{z_0}}\ln\left|\frac{\sqrt{z_2}-\sqrt{z_0}}{\sqrt{z_2}+\sqrt{z_0}}\frac{\sqrt{z_1}+\sqrt{z_0}}{\sqrt{z_1}-\sqrt{z_0}}\right|.\label{ap1d2}
\end{align}

\item Finally, when $n=1$ the elements are
\begin{align}
M_{1,0}=&-\frac{2}{\sqrt{z_1}+\sqrt{z_0}},\label{ap1e1}\\
M_{1,1}=&\frac{4\tau}{(\sqrt{z_1}+\sqrt{z_0})(\sqrt{z_2}+\sqrt{z_0})(\sqrt{z_2}+\sqrt{z_1})}-\nonumber\\
&\frac{1}{\sqrt{z_1}}\ln\left|\frac{\sqrt{z_2}+\sqrt{z_1}}{\sqrt{z_2}-\sqrt{z_1}}\frac{\sqrt{z_0}-\sqrt{z_1}}{\sqrt{z_0}+\sqrt{z_1}}\right|.\label{ap1e2}
\end{align}
\item Again in case $n=0,1$ and $m\ne 0,1$ the elements are determined by general formulas (\ref{ap1a1})--(\ref{ap1a7}). 
\end{itemize}
The elements $(0,0)$ and $(N,N)$ are formally infinite diverging logarithmically but are assumed to be finite and equal to $-10^8$ and $10^8$, respectively, as it does not substantially change the results of inversion.

Finally, the elements of vectors $\overrightarrow{\mathbf{g}}$ are calculated from (\ref{2k2}) using the trapezoidal rule 
\begin{multline}
g_n=\sqrt{\frac{2ki}{\pi z_n}}\sum\limits_{s=1}^{N_x-1}u(x_s,z_n)\exp\left[-\frac{ik{x_s}^2}{2z_n}\right]+\nonumber\\
\frac{1}{2}\sqrt{\frac{2ki}{\pi z_n}}\left(u(x_0,z_n)\exp\left[-\frac{ik{x_0}^2}{2z_n}\right]+u(x_{N_x},z_n)\exp\left[-\frac{ik{x_{N_x}}^2}{2z_n}\right]\right)\label{ap1g1}
\end{multline}
where $N_x+1$ is the number of points in the image line where the amplitude is measured, $x_s=x_{min}+sh$, $h=(x_{max}-x_{min})/N_x$ and $[x_{min},x_{max}]$ is the interval in the image line where the amplitude is known.

\end{document}